\documentclass[10pt,twoside]{article}
\usepackage{amssymb}
\usepackage{mathrsfs}
\usepackage{amsmath}
\usepackage{latexsym}
\newtheorem{theorem}{Theorem}[section]

\newtheorem{lemma}{Lemma}[section]
\newtheorem{proposition}{Proposition}[section]
\newtheorem{remark}{Remark}[section]

\newtheorem{corollary}{Corollary}[section]
\def\proof{\mbox {\it Proof.~}}

\newcommand{\N}{\mathbb{N}}                     
\newcommand{\R}{\mathbb{R}}                     
\newcommand{\set}[2]{\left\{{#1}\mid{#2}\right\}}       
\newcommand{\qed}{\hfill $\Box$ \bigskip}       
\newcommand{\dist}{\mathrm{dist\,}}             
\newcommand{\ind}{\mathrm{ind\,}}               
\newcommand{\crit}{\mathrm{crit}\,}               
\newcommand{\grad}{\mathrm{grad\,}}             
\newcommand{\sing}{\mathrm{sing}\,}             

\makeatletter\def\theequation{\arabic{section}.\arabic{equation}}\makeatother

\begin{document}
\title{
\vspace{0.5in} {\bf\Large  A Smooth Pseudo-Gradient for the
Lagrangian Action Functional}}

\author{
{\bf\large Alberto Abbondandolo}\footnote{The first author was partially supported by a {\em Humboldt Research Fellowship for Experienced Researchers}.}\hspace{2mm}
{\bf\large }\vspace{1mm}\\
{\it\small Dipartimento di Matematica}\\ {\it\small Universit\`a di
Pisa},
{\it\small I-56127 Pisa, Italy}\\
{\it\small e-mail: abbondandolo@dm.unipi.it\vspace{2mm}}\\
{\bf\large Matthias Schwarz}\footnote{The second author was
partially supported by the DFG grant SCHW 892/2-3.}\hspace{2mm}\\
\vspace{1mm} {\it\small Mathematisches Institut}\\ {\it\small
Universit\"at Leipzig},
{\it\small D-04103 Leipzig, Germany}\\
{\it\small e-mail: mschwarz@mathematik.uni-leipzig.de\vspace{2mm}}\\
{\small\bf To Vieri, with admiration and gratitude}}

\date{November 4, 2009}

\maketitle

\begin{center}
{\bf\small Abstract}

\vspace{3mm}
\hspace{.05in}\parbox{4.5in}
{{\small We study the action functional associated to a smooth Lagrangian function on the tangent bundle of a manifold, having quadratic growth in the velocities. We show that, although the action  functional is in general not twice differentiable on the Hilbert manifold consisting of $H^1$ curves, it is a Lyapunov function for some smooth Morse-Smale vector field, under the generic
assumption that all the critical points are non-degenerate.
This fact is sufficient to associate a Morse complex to the Lagrangian
action functional.
}}
\end{center}


\renewcommand{\theenumi}{\roman{enumi}}
\renewcommand{\labelenumi}{(\theenumi)}

\section{\bf Introduction}
\def\theequation{1.\arabic{equation}}\makeatother
\setcounter{equation}{0}

This paper is about the action functional
\[
\mathbb{S}_L(\gamma) = \int_0^1 L(t,\gamma(t),\gamma'(t))\, dt
\]
associated to a smooth Lagrangian function $L$ on $[0,1]\times TM$, where $M$ is a smooth $n$-dimensional manifold, and $\gamma$ is a curve in $M$. The
critical points of such a functional are the solutions of the Euler-Lagrange equation associated to $L$.
If $L$ has quadratic growth in the velocities (see condition (L1) in Section \ref{stabe}), the natural functional setting for studying the functional $\mathbb{S}_L$ is the Hilbert manifold $H^1([0,1],M)$ of absolutely continuous paths in $M$ with square integrable derivative, or the Hilbert submanifolds one gets by imposing the boundary conditions one is interested in. In fact, as Benci has showed in \cite{ben86}, this assumption implies that the action functional $\mathbb{S}_L$ is continuously differentiable on $H^1([0,1],M)$ and, assuming also that $L$ is fiber-wise uniformly convex and $M$ is compact, that the Palais-Smale condition holds (see Section \ref{stabe} for precise statements). These facts can be used to prove - for instance - that if $M$ is compact the Euler-Lagrange equation associated to a time-periodic $L$ has at least as many periodic solutions as the Lusternik-Schnirelmann category of the free loop space of $M$ (see \cite{ben86}).

Under this quadratic growth assumption on $L$, the functional $\mathbb{S}_L$ is actually also twice Gateaux differentiable\footnote{In \cite{af07} and \cite{as06} it is erroneously stated that the action functional $\mathbb{S}_L$ is indeed $C^2$. Actually, $C^1$ regularity is sufficient in \cite{af07}, which deals with Lusternik-Schnirelmann theory. The error in \cite{as06} was corrected in \cite{aps08}, by considering a smaller class of Lagrangian functions. The main result of this paper implies that the original class of Lagrangian functions can also be considered in \cite{as06}.} on $H^1([0,1],M)$, but it is twice
Fr\'ech\'et differentiable if and only if $L$ is {\em exactly} quadratic in the velocities (more precisely, if for every $(t,q)$ the function $v\mapsto L(t,q,v)$ is a polynomial of degree at most 2). In this case, $\mathbb{S}_L$ is actually smooth. This fact is related to well-known facts about the differentiability of Nemitsky operators. Infinite dimensional Morse theory for the functional $\mathbb{S}_L$ (see e.g. \cite{cha93}), and in particular the Morse lemma, would require $C^2$ regularity.

The aim of this paper is to show that if the critical points of $\mathbb{S}_L$ are non-degenerate, even if $\mathbb{S}_L$ is not $C^2$, there exists a {\em smooth} (i.e.\ $C^{\infty}$) Morse-Smale vector field $X$ on $H^1([0,1],M)$ for which $\mathbb{S}_L$ is a Lyapunov function. The existence of such a vector field allows to construct the {\em Morse complex} of $\mathbb{S}_L$, a chain complex on the free Abelian group generated by the critical points of $\mathbb{S}_L$, graded by their Morse index. The homology of such  a chain complex is isomorphic to the singular homology of the underlying path space.

These facts imply multiplicity results for the solution of the Euler-Lagrange equation associated to $L$, under various boundary conditions (for instance, Dirichlet, Neumann, or periodic boundary conditions). Here we treat these various boundary conditions in a unified way, by fixing a submanifold $Q$ of $M\times M$ and by restricting the action functional $\mathbb{S}_L$ to the submanifold $H^1_Q([0,1],M)$ of $H^1([0,1],M)$ consisting of those curves $\gamma$ such that $(\gamma(0),\gamma(1))$ belongs to $Q$. See Corollary
\ref{coro} below.

These existence results could be proved also by determining the
critical groups of the critical points of $\mathbb{S}_L$ by some
weak version of the Morse Lemma (see \cite{lll05} and \cite{dhk07}
for two abstract statements of this kind, and \cite{cjm09} for their
use in the study of the Finsler energy functional), or by replacing
the space $H^1([0,1],M)$ by finite dimensional spaces of continuous
piecewise extremals of the Euler-Lagrange equation (see \cite{rad92}
for the case where $L$ is the Finsler energy, or \cite{mazz08} for
more general Lagrangian functions). However, our main motivation for
constructing the Morse complex for the action functional associated
to a Lagrangian $L$ which has quadratic growth in the velocities
comes from the study of Floer homology of cotangent bundles and its
relationship with string topology. As shown in \cite{as06} and
\cite{aps08}, the homology of the Floer complex associated to a
Hamiltonian vector field on the cotangent bundle of a closed
manifold and to the non-local conormal boundary condition induced by
a submanifold $Q$ of $M\times M$, can be determined by coupling the
Cauchy-Riemann type equation arising in Floer theory with the
negative gradient flow of $\mathbb{S}_L$ on $H^1_Q([0,1],M)$. In
some situations, one can limit the attention to Lagrangian functions
$L$ which are exactly quadratic in the velocities, so that the
gradient of $\mathbb{S}_L$ is indeed smooth. However, in other cases
it is useful to consider more general Lagrangian functions, for
instance because this allows to achieve the transversality needed in
Floer theory without perturbing the Levi-Civita almost complex
structure on $T^*M$, as in \cite{as08}. With such a Lagrangian
function $L$, the negative gradient flow of $\mathbb{S}_L$ is not
regular enough (the stable and unstable manifolds of critical points
are just Lipschitz submanifolds), but the smooth vector field that
we construct in this paper turns out to be a good replacement for the
negative gradient vector field.\\
\indent This paper is organized in the following way. In Section
\ref{morcom} we recall the construction of the Morse complex
associated to a $C^1$ function on a Hilbert manifold, which is a
Lyapunov function for some Morse-Smale $C^1$ vector field having
singular points with finite Morse index. In Section \ref{stabe} we study
the differentiability and compactness properties of the Lagrangian
action functional on $H^1([0,1],M)$. Most of the results of this
section are well-known, but for sake of completeness we include
complete proofs. In particular, we prove the not so well-known fact
that if $\mathbb{S}_L$ is twice Fr\'ech\'et differentiable at some
curve $\gamma\in H^1([0,1],M)$, then for every $t$ the function
$v\mapsto L(t,\gamma(t),v)$ is a polynomial of degree at most 2 on
$T_{\gamma(t)} M$. In Section 4 we prove the existence of a smooth
Morse-Smale vector field $X$ for which $\mathbb{S}_L$ is a Lyapunov
function. The main point is to construct $X$ in a neighborhood of
each critical point. This is done in Lemma \ref{loc}, by using some
properties of the Gateaux second differential of $\mathbb{S}_L$,
together with the strong continuity of the functional calculus on
the space of self-adjoint bounded operators on a Hilbert space, and
other facts about the strong convergence of self-adjoint bounded
operators.

\paragraph{Acknowledgment.} The first author wishes to thank the Max-Planck-Institut f\"ur Mathematik in den Naturwissenschaften of Leipzig for its kind hospitality.
\newpage
\section{\bf The Morse complex}
\def\theequation{2.\arabic{equation}}\makeatother
\setcounter{equation}{0}

\label{morcom}

The aim of this section is to recall the construction and the basic properties of the Morse complex. Detailed proofs can be found in \cite{ama06m} (in a more general setting).

Let $\mathscr{M}$ be a smooth Hilbert manifold. This means that $\mathscr{M}$ is a paracompact Hausdorff topological space which is locally homeomorphic to a real Hilbert space $\mathbb{H}$, and which admits an atlas whose
transition maps are smooth (see \cite{lan99} for foundational issues concerning infinite dimensional manifolds). Let $X$ be a $C^k$ vector field on $\mathscr{M}$, with $1\leq k \leq \infty$. The associated local flow is the $C^k$ map $(t,p) \mapsto \phi(t,p)$ obtained by solving the Cauchy problem
\[
\frac{\partial}{\partial t} \phi(t,p) = X(\phi(t,p)), \quad \phi(0,p) = p,
\]
which is locally well-posed because $X$ is at least $C^1$.
The vector field $X$ is said to be {\em complete} (respectively {\em forward complete} or {\em backward complete}) if the maximal domain of $\phi$ is the whole $\R \times \mathscr{M}$ (respectively, if it contains $[0,+\infty[ \times \mathscr{M}$ or $]-\infty,0] \times \mathscr{M}$).

The set of singular points of $X$, which is denoted by $\sing(X)$, is the set of points $x\in \mathscr{M}$ such that $X(x)=0$. If $x\in \sing(X)$, the differential of $X$ at $x$, which is a continuous linear map from $T_x \mathscr{M}$ to $T_{X(x)} T \mathscr{M}$, induces a bounded operator on $T_x \mathscr{M}$, that we denote by $\nabla X(x)$. The singular point $x$ is said to be {\em hyperbolic} if the spectrum of $\nabla X(x)$ is disjoint from the imaginary axis $i\R$. This implies, in particular, that $\nabla X(x)$ is invertible, so hyperbolic singular points are isolated points of $\sing(X)$.
A vector field all of whose singular points are hyperbolic is said to be {\em Morse}.
The division of the spectrum of $\nabla X(x)$ into the component with negative real part and the one with positive real part induces a $\nabla X(x)$-invariant splitting of $T_x \mathscr{M}$ into two closed linear subspaces, that we denote by
\[
T_x \mathscr{M} = E^s(x) \oplus E^u(x).
\]
The {\em Morse index} $\ind(x)$ of $x$ is the (possibly infinite) dimension of $E^u(x)$, the invariant subspace corresponding to the positive part of the spectrum. The stable manifold theorem implies that the subsets
\begin{eqnarray*}
W^s(x) & := & \set{p\in \mathscr{M}}{\phi(t,p) \rightarrow x \mbox{ for } t\rightarrow +\infty}, \\
W^u(x) & := & \set{p\in \mathscr{M}}{\phi(t,p) \rightarrow x \mbox{ for } t\rightarrow -\infty}
\end{eqnarray*}
are immersed $C^k$ submanifolds of $\mathscr{M}$. They are called the {\em stable} and {\em unstable manifold of $x$}, respectively. Actually, they are images of $E^s(x)$ and $E^u(x)$ under injective $C^k$ immersions which map $0$ to $x$ and whose differential at $0$ is the inclusion mapping (see \cite{ama06} for the question of the global differentiable structure of the stable and unstable manifold, in a Banach space setting).

A $C^1$ function $f:\mathscr{M} \rightarrow \R$ is said to be a Lyapunov function for $X$ if
\[
Df(p)[X(p)] < 0 , \quad \forall p\in \mathscr{M}\setminus \sing(X).
\]
The set of critical points of the Lyapunov function $f$, that we denote by $\crit(f)$, is contained in $\sing(X)$, and if $X$ is Morse then $\crit (f) = \sing(X)$. If the hyperbolic singular point $x$ has finite Morse index and if $X$ admits a Lyapunov function, then the stable and unstable manifolds $W^s(x)$ and $W^u(x)$ are locally closed $C^k$ {\em embedded} submanifolds (if we drop the assumption on the finiteness of the Morse index, the same fact is true under a mild non-degeneracy condition on $f$, see Theorem 1.20 in \cite{ama06m}).

Assume that $X$ is Morse, that all its singular points have finite Morse index, and that $X$ admits a Lyapunov function $f$. Then $X$ is said to satisfy the {\em Morse-Smale property up to order $h$} if the intersection of $W^u(x)$ with $W^s(y)$ is transverse, whenever $x$ and $y$ are singular points with $\ind(x)-\ind(y) \leq h$. In particular, the Morse-Smale property up to order 0 amounts to asking that $W^u(x)\cap W^s(y)$ should be empty whenever $\ind(x)\leq \ind(y)$ and $x\neq y$.

A sequence $(p_h)\subset \mathscr{M}$ is said to be a {\em Palais-Smale sequence} for the pair $(X,f)$ if $f(p_n)$ is bounded and $Df(p_n)[X(p_n)]$ is infinitesimal. The pair $(X,f)$ is said to satisfy the {\em Palais-Smale condition} if every Palais-Smale sequence has a converging subsequence. When $\mathscr{M}$ is endowed with a Riemannian metric and $-X$ is the gradient of $f$ with respect to such a metric, we recover the standard definitions (see e.g.\ \cite{cha93}).

Let us assume that the vector field $X$ satisfies the following conditions:
\begin{enumerate}
\item[(X1)] $X$ is of class $C^1$;
\item[(X2)] $X$ is Morse and all its singular points have finite Morse index;
\item[(X3)] $X$ admits a $C^1$ Lyapunov function $f$, which is bounded from below;
\item[(X4)] the pair $(X,f)$ satisfies the Palais-Smale condition;
\item[(X5)] $X$ is forward complete;
\item[(X6)] $X$ satisfies the Morse-Smale condition up to order zero.
\end{enumerate}

Let us fix a real number $a$.
Since $f$ is bounded from below, assumptions (X2) and (X4) imply that $X$ has finitely many singular points in the sublevel $\{f<a\}$. Every singular point $x\in \{f<a\}$ admits a suitable small open neighborhood $\mathscr{U}_x$ such that, setting
\begin{equation}
\label{cellfil}
\mathscr{M}_k^a := \bigcup_{\substack{x\in \sing(X)\\ f(x)< a \\ \ind(x) \leq k}} \phi([0,+\infty[ \times \mathscr{U}_x),
\quad \forall k\in \N \cup \{\infty\},
\end{equation}
the sequence of open sets $\{\mathscr{M}_k^a\}_{k\in \N}$ is a cellular filtration of $\mathscr{M}_{\infty}^a$ (such a filtration is eventually constant, since $\{f<a\}$ contains finitely many singular points). More precisely, denoting by $H_*$ the singular homology groups with integer coefficients, we have
\begin{equation}
\label{ceco}
H_j(\mathscr{M}_k^a,\mathscr{M}_{k-1}^a) \cong \left\{ \begin{array}{ll} M_k^a(X) & \mbox{if } j=k, \\ 0 & \mbox{if } j\neq k, \end{array} \right.
\end{equation}
where $M_k^a(X)$ is the free Abelian group generated by the singular points $x$ of $X$ with Morse index $\ind(x)=k$ and $f(x)<a$. In fact, $H_*(\mathscr{M}_k^a,\mathscr{M}_{k-1}^a)$ is generated by the $k$-dimensional relative homology classes induced by the embeddings of the local unstable manifolds $W^u(x)\cap \mathscr{U}_x$ by the map $\phi(t,\cdot)$, for $t>0$ large enough, for any singular point $x\in \{f<a\}$ of Morse index $k$. Therefore, the isomorphism $H_k(\mathscr{M}_k^a,\mathscr{M}_{k-1}^a) \cong M_k^a(X)$ is induced by the choice of an orientation of each unstable manifold $W^u(x)$, for $x\in \sing(X) \cap \{f<a\}$ of Morse index $k$.
Furthermore, the open set $\mathscr{M}_{\infty}^a$ is a deformation retract of the sublevel $\{ f< a\}$.
By (\ref{ceco}), the chain complex associated to the cellular filtration $\{\mathscr{M}_k^a \}_{k\in \N}$ consists of the graded group $M_*^a(X)$ with the boundary homomorphism
\[
\partial_k : M_k^a(X) = H_k(\mathscr{M}_k^a,\mathscr{M}_{k-1}^a) \rightarrow H_{k-1}(\mathscr{M}_{k-1}^a,\mathscr{M}_{k-2}^a) = M_{k-1}^a(X).
\]
The chain complex $\{M_*^a(X),\partial_*\}$ is called the {\em Morse complex of $X$ on the sublevel $\{f<a\}$}. The boundary homomorphism $\partial_*$ does not change if we choose different neighborhoods $\mathscr{U}_x$.
By the standard isomorphism between cellular homology and singular homology (see e.g.\ \cite{dol80}), the homology of the Morse complex of $X$ on $\{f<a\}$ is isomorphic to the singular homology of $\mathscr{M}_{\infty}^a$, hence to the the singular homology of the sublevel $\{f<a\}$.

We have so far restricted the attention to a  sublevel $\{f<a\}$, because in general one may not be able to find open neighborhoods $\mathscr{U}_x$ of the critical points $x$ in such a way that the sets defined as in (\ref{cellfil}) but without the upper bound on $f$ constitute a cellular filtration. The reason is that, in general, (X6) does not exclude that a singular point $x$ may belong to the closure of the union of all singular points with Morse index not larger than $\ind(x)$ (unless $X$ has finitely many singular points of any given Morse index). However, the {\em Morse complex of $X$}, $\{M_*(X),\partial_*\}$, can be defined as the direct limit of the chain complexes $\{M_*^a(X),\partial_*\}$, for $a\uparrow +\infty$. Since the homology of a direct limit of complexes is the
direct limit of the homologies, the homology of $\{M_*(X),\partial_*\}$ is
isomorphic to the singular homology of $\mathscr{M}$.

Let us strengthen condition (X6) by assuming that:
\begin{enumerate}
\item[(X6')] $X$ satisfies the Morse-Smale condition up to order one.
\end{enumerate}
Then the boundary operator $\partial_*$ can be expressed in terms of suitable intersection numbers. The chosen orientation of each unstable manifold $W^u(x)$ determines a co-orientation - that is, an orientation of the normal bundle - of each stable manifold $W^s(x)$. Since the transverse intersection of an oriented and a co-oriented submanifold carries a canonical orientation, by (X6') we get an orientation for each intersection $W^u(x)\cap W^s(y)$, when $\ind(x)-\ind(y)=1$. In this case, $W^u(x)\cap W^s(y)$ consists of finitely many flow lines, and the integer $n(x,y)$ can be defined as the number of those lines on which this orientation agrees with the direction of the flow minus the number of the other lines. The boundary homomorphism $\partial_*$ can be expressed in terms of the generators $x\in \sing(X)$ of $M_*(X)$ as
\begin{equation}
\label{form}
\partial_k x = \sum_{\substack{y\in \sing(X)\\ \ind(y)=k-1}} n(x,y)\, y,
\end{equation}
for every $x\in \sing(X)$ with $\ind(x)=k$.

If two vector fields $X_1$ and $X_2$ satisfy (X1)-(X6) with respect
to the same Lyapunov function $f$, then they have the same singular points,
with the same Morse indexes, and there is a natural isomorphism
between their Morse complexes. In general, this isomorphism will not
be  the identity mapping on $M_*(X_1)=M_*(X_2)$, which need not be a
chain map, but it can be defined by applying the flow of the vector
field $X_2$ to the cellular filtrations associated to $X_1$.\\
\indent Let $f$ be a $C^1$ function on $\mathscr{M}$ which admits a
vector field $X$ satisfying (X1)-(X5) (with respect to $f$ itself).
Then $X$ can be perturbed in such a way that also (X6) holds (and
even (X6')). These Morse-Smale perturbations are generic in suitable
complete metric spaces consisting of $C^1$ vector fields $Y$
agreeing with $X$ at each singular point up to order one, and endowed
with a Whitney metric involving the function $|Df[Y]|$, such that
any $Y$ whose distance from $X$ is smaller than 1 has the same rest
points of $X$ and satisfies (X1)-(X5).\\
\indent These facts imply that every $C^1$ function $f$ which admits
a vector field $X$ satisfying (X1)-(X5) has a Morse complex, which
is uniquely defined up to isomorphism. If $f$ is a $C^2$ function on
$\mathscr{M}$ which is bounded from below, has only non-degenerate
critical points with finite Morse index, and satisfies the
Palais-Smale condition  (in the usual sense) with respect to a
complete Riemannian metric on $\mathscr{M}$, then the negative
gradient of $f$ satisfies (X1)-(X5). In particular, such functions
have a Morse complex, defined up to isomorphism. As this paper shows,
there are interesting examples of functions which are $C^1$ and not
$C^2$, but for which there exists a smooth vector field which
satisfies (X1)-(X5) and the Morse-Smale condition up to every order.
\vspace{3mm}
\begin{remark}
\label{morreg} {\rm The proof that a vector field $X$ which
satisfies (X1)-(X5) can be perturbed in a generic way in order to
achieve also (X6') is based on the Sard-Smale theorem. The
Sard-Smale theorem for a Fredholm map $\Phi$ requires $\Phi$ to be
$C^k$, with $k$ at least one and larger than the Fredholm index of
$\Phi$. Therefore, by perturbing in a generic way a vector field $X$
which is just $C^1$, one obtains only the Morse-Smale condition only
up to order 1 (by using the $\R$-action, the Fredholm map one needs
to consider in this case has index zero). More generally, a generic
perturbation in the $C^k$ topology of a $C^k$ Morse vector field
satisfies the Morse-Smale condition up to order $k$, for every $k\in
\N \cup\{\infty\}$. We also recall that if $k\geq 1$, $C^{k+1}$
functions on an infinite dimensional Hilbert space are not dense in
the space of $C^k$ functions (see \cite{ll86}). These considerations
are relevant only when the object of interest is the vector field.
When we start from a $C^1$ function $f$ admitting a $C^1$ vector
field $X$ satisfying (X1)-(X5), it is often easy to approximate $X$
in the $C^0$ topology by a $C^{\infty}$ vector field which satisfies
the same conditions, which can then be further perturbed so to
achieve the Morse-Smale condition up to every order.}
\end{remark}
\newpage
\section{\bf Basic properties of the Lagrangian action \newline functional}
\def\theequation{3.\arabic{equation}}\makeatother
\setcounter{equation}{0}

\label{stabe}

In this section we recall some well known facts about Hilbert
manifolds of paths and about the properties of the Lagrangian
action functional on such Hilbert manifolds. A standard reference
for the former issue is \cite{kli82}, while the latter one is
treated in \cite{ben86}. For sake of completeness, we include
detailed proofs, simplifying some of the arguments in the above
mentioned references.

\paragraph{Lagrangian functions with quadratic growth.} Let $(M,g)$ be a (not necessarily compact or complete) smooth $n$-dimensional Riemannian manifold without boundary.
Points in $TM$ are denoted by $(q,v)$, where $q\in M$ and $v\in T_q M$. The metric $g$ induces a metric on $TM$, Levi-Civita covariant derivations both on $M$ and on $TM$, the horizontal subbundle  $T^h TM$ of $TTM$, and isomorphisms
\[
T_{(q,v)} TM = T^h_{(q,v)} TM \oplus  T^v_{(q,v)} TM \cong T_q M \oplus T_q M,
\]
where $ T^v_{(q,v)} TM = \ker D\pi(q,v)$ is the vertical subspace, $\pi: TM \rightarrow M$ being the canonical projection. The horizontal and vertical component of the gradient of a function defined on $TM$ are denoted by $\nabla_q$ and $\nabla_v$, respectively. We use a similar notation for higher order derivatives.

The Lagrangian functions we consider are smooth functions on $[0,1]\times TM$ which satisfy the following conditions:
\begin{enumerate}
\item[(L1)] There is a continuous function $\ell_1$ on $M$ such that
\begin{eqnarray*}
\left\| \nabla_{vv} L (t,q,v) \right\| & \leq & \ell_1(q), \\ \left\|
\nabla_{vq} L (t,q,v) \right\|&  \leq & \ell_1(q) (1+g(v,v)^{1/2}),  \\
\left\| \nabla_{qq} L (t,q,v) \right\| & \leq & \ell_1(q) (1+g(v,v)),
\end{eqnarray*}
for every $(t,q,v)\in [0,1]\times TM$.
\item[(L2)] There is a continuous positive function $\ell_2$ on $M$ such that $\nabla_{vv} L (t,q,v) \geq \ell_2(q) I$, for every $(t,q,v)\in [0,1]\times TM$.
\end{enumerate}
Condition (L1) implies that $L$ grows at most quadratically on each fiber. Condition (L2) implies that $L$ is fiber-wise strictly convex and grows at least quadratically on each fiber. These conditions are independent on the choice of the metric $g$, in the sense that if $L$ satisfies (L1) and (L2) with respect to the metric $g$, then it satisfies the same conditions with respect to any other metric, with different functions $\ell_1$ and $\ell_2$. They are semi-local conditions, in the sense that they are local in $q\in M$ but they involve all values of $v\in T_q M$.
Classical examples of Lagrangian functions which satisfy (L1) and (L2) are the so called {\em electro-magnetic} Lagrangian functions,
that is functions of the form
\begin{equation}
\label{em}
L(t,q,v) = \frac{1}{2} \langle A(t,q) v,v \rangle + \langle \alpha(t,q), v \rangle - V(t,q),
\end{equation}
where $\langle\cdot,\cdot \rangle$ denotes the duality paring, $A(t,q) : T_q M \rightarrow T^*_q M$ is a symmetric, strictly positive linear mapping, smoothly depending on $(t,q)$ (the kinetic tensor), $\alpha$ is a smoothly $t$-dependent one-form (the magnetic potential), and $V$ is a smooth function on $[0,1]\times M$ (the potential energy).

It is easy to translate the assumptions (L1) and (L2) in terms of local coordinates. Actually, it is useful to work with {\em time-dependent} local coordinates systems, that is smooth maps
\[
\varphi: [0,1] \times U \rightarrow M,
\]
where $U$ is an open subset of $\R^n$, and $\varphi(t,\cdot)$ is a diffeomorphism onto the open subset $\varphi(\{t\} \times U)$, for every $t\in [0,1]$. We denote by $\mathscr{F}$ the family of all such pairs $(U,\varphi)$.
It is often useful to assume the element $(U,\varphi)$ to be {\em bi-bounded}, meaning that $U$ is bounded, $\varphi([0,1]\times U)$ has compact closure in $M$, and all the derivatives of $\varphi$ and of the map $(t,q) \mapsto \varphi(t,\cdot)^{-1}(q)$ are bounded. Every element $(U,\varphi)$ induces the time-dependent coordinate system
\[
[0,1]\times U \times \R^n \rightarrow TM, \quad (t,q,v) \mapsto (\varphi(t,q),D_q \varphi(t,q)[v]),
\]
on $TM$. The pull-back of $L$ by such a coordinate system is the function
\[
(\varphi^* L) (t,q,v) = L(t,\varphi(t,q), D_q \varphi(t,q)[v]), \quad \forall (t,q,v) \in [0,1] \times U \times \R^n.
\]
When no confusion is possible, we denote $\varphi^* L$ simply by $L$. With this convention, conditions (L1) and (L2) can be restated by saying that for every bi-bounded coordinate system $(U,\varphi)$ there hold:
\begin{enumerate}
\item[(L1')] There is a number $\ell_1$ such that
\[
\hspace{-.8cm}
{\textstyle
\left\| \frac{\partial^2 L}{\partial v^2} (t,q,v) \right\|  \leq  \ell_1, \quad \left\|
\frac{\partial^2 L}{\partial v \partial q} (t,q,v) \right\|  \leq  \ell_1 (1+ |v|),  \quad
\left\| \frac{\partial^2 L}{\partial q^2} (t,q,v) \right\|  \leq  \ell_1 (1+ |v|^2),}
\]
for every $(t,q,v)\in [0,1]\times U \times \R^n$.
\item[(L2')] There is a positive number $\ell_2$ such that $\partial^2 L/\partial v^2 (t,q,v) \geq \ell_2 I$, for every $(t,q,v)\in [0,1]\times U \times \R^n$.
\end{enumerate}
If integrated along the fibers, (L1') implies the growth conditions
\begin{eqnarray}
\label{Agc1}
\left| \frac{\partial L}{\partial v} (t,q,v) \right| \leq \ell_3 (1 + |v|), \quad
\left| \frac{\partial L}{\partial q} (t,q,v) \right| \leq \ell_3 (1+|v|^2),
\\ \label{Agc2}
L(t,q,v) \leq \ell_4 (1 + |v|^2),
\end{eqnarray}
for suitable numbers $\ell_3$ and $\ell_4$.

\paragraph{The non-local boundary value problem.} By (L2), the Euler-Lagrange equation associated to the Lagrangian function $L$, which in local coordinates can be written as
\begin{equation}
\label{Aele}
\frac{d}{dt} \left( \frac{\partial L}{\partial v} (t,\gamma(t),\gamma'(t)) \right) = \frac{\partial L}{\partial q} (t,\gamma(t),\gamma'(t)),
\end{equation}
defines a locally well-posed second order Cauchy problem. We shall treat different boundary conditions for equation (\ref{Aele}) in a unified way, by considering a smooth submanifold (non-empty, without boundary, not necessarily closed) $Q$ of $M\times M$, and by imposing the non-local boundary conditions
\begin{eqnarray}
\label{Abdry1}
 (\gamma(0),\gamma(1)) \in Q,  \\ \label{Abdry2}
D_v L(0,\gamma(0),\gamma'(0)) [\xi_0] = D_v L(1,\gamma(1),\gamma'(1)) [\xi_1], \quad \forall (\xi_0,\xi_1) \in T_{(\gamma(0),\gamma(1))} Q.
\end{eqnarray}
Here $D_v L(t,q,v) : T_q M \rightarrow \R$ is the fiber-wise
differential. \vspace{2mm}
\begin{remark}
{\rm If $H:[0,1]\times T^*M\rightarrow \R$ is the Fenchel-dual of
$L$, the boundary value problem
(\ref{Aele}-\ref{Abdry1}-\ref{Abdry2}) is equivalent to the problem
of finding Hamiltonian orbits $x:[0,1]\rightarrow T^*M$ such that
the pair $(x(0),-x(1)) \in T^*M \times T^*M = T^* (M\times M)$
belongs to the conormal bundle of $Q$, that is to the set of all
covectors in $T^* (M\times M)$ which are based at points of $Q$ and
annihilate the tangent space of $Q$.}
\end{remark}
\paragraph{Hilbert manifolds of Sobolev paths.}
Let us consider the set $H^1([0,1],M)$ consisting of all absolutely continuous curves $\gamma:[0,1]\rightarrow M$ with square-integrable first derivative. It is well-known that this set has a natural structure of a smooth Hilbert manifold. In order to study local properties of the Lagrangian action functional, it is convenient to recall the construction of this structure here.

Let $(U,\varphi)$ be an element of $\mathscr{F}$.
The set $H^1([0,1],U)$, consisting of all the elements $\gamma$ of the Hilbert space $H^1([0,1],\R^n)$ whose image lies in $U$, is open, because of the continuity of the embedding $H^1 \hookrightarrow C^0$. The map $\varphi$ induces an injective map
\[
\varphi_* : H^1([0,1],U) \rightarrow H^1([0,1],M), \quad \varphi_*(\gamma) := \varphi(\cdot,\gamma(\cdot)),
\]
whose image is denoted by $\mathscr{U}_{(U,\varphi)}$. The Hilbert manifold structure on $H^1([0,1],M)$ is defined by declaring the family of maps $\varphi_*$, for $(U,\varphi)$ in $\mathscr{F}$, to be an atlas. The sets $\mathscr{U}_{(U,\varphi)}$ are actually open also with respect to the $C^0$ topology on $H^1([0,1],M)$. Notice that every element $\gamma$ of $H^1([0,1],M)$ belongs to $\mathscr{U}_{(U,\varphi)}$ for some $(U,\varphi)$ in $\mathscr{F}$. Indeed, if $\tilde{\gamma}$ is a smooth curve in $M$ such that $\dist(\gamma(t),\tilde{\gamma}(t))$ is small enough, then we can define a map $\varphi$ as above and such that $\gamma\in \mathscr{U}_{(U,\varphi)}$ by composing the map
\[
(t,\xi) \mapsto \exp_{\tilde{\gamma}(t)} (\xi), \quad \xi \in T_{\tilde{\gamma}(t)} M,
\]
with a smooth trivialization of the vector bundle $\tilde{\gamma}^*(TM)$ over $[0,1]$.

The subset
\[
H^1_Q([0,1],M) = \set{\gamma \in H^1([0,1],M)}{(\gamma(0),\gamma(1)) \in Q}
\]
is a smooth submanifold, being the inverse image of $Q$ by the smooth submersion
\[
H^1([0,1],M) \rightarrow M \times M, \quad \gamma \mapsto (\gamma(0),\gamma(1)).
\]
Actually, a smooth atlas for $H^1_Q([0,1],M)$ can be build by fixing a linear subspace $W$ of $\R^n \times \R^n$ with $\dim W = \dim Q$, by considering elements $(U,\varphi)\in \mathscr{F}$ such that $0\in U$ and $(\varphi(0,q),\varphi(1,q))$ belongs to $Q$ for every $q$ in $U\cap W$, and by restricting the map $\varphi_*$ to the intersection of the open set $H^1([0,1],U)$ with the closed linear subspace
\[
H^1_W([0,1],\R^n) := \set{ \gamma \in H^1([0,1],\R^n)}{(\gamma(0),\gamma(1))\in W}.
\]
The set of pairs $(U,\varphi)\in \mathscr{F}$ with the above property is denoted by $\mathscr{F}_Q$.

\paragraph{Differentiability of the action functional.} If $L$ satisfies (L1), the Lagrangian action functional
\[
\mathbb{S}_L(\gamma) := \int_0^1 L(t,\gamma(t),\gamma'(t))\, dt
\]
is well-defined on $H^1([0,1],M)$, because of (\ref{Agc2}).  The restriction of the functional $\mathbb{S}_L$ to the submanifold $H^1_Q([0,1],M)$ is denoted by $\mathbb{S}_L^Q$.

If $(\varphi,U)$ belongs to $\mathscr{F}$, then
\[
\mathbb{S}_L(\varphi_*(\gamma)) = \mathbb{S}_{\varphi^*(L)} (\gamma), \quad \forall \gamma \in H^1([0,1],U),
\]
so the study of the local properties of $\mathbb{S}_L$ is reduced to the study of the functional $\mathbb{S}_{\varphi^*(L)}$, which is defined on an open subset of a Hilbert space.

\begin{proposition}
\label{reg}
Assume that the smooth Lagrangian function $L:[0,1]\times TM \rightarrow \R$ satisfies (L1). Then:

\begin{enumerate}

\item The Lagrangian action functional $\mathbb{S}_L$ is continuously differentiable on the Hilbert manifold $H^1([0,1],M)$. Its differential $D\mathbb{S}_L$ is locally Lipschitz continuous and Gateaux-differentiable.

\end{enumerate}

\noindent
 Let $Q$ be a submanifold of $M\times M$, and assume that $L$ satisfies also (L2). Then:

\begin{enumerate}

\setcounter{enumi}{1}

\item The critical points of $\mathbb{S}_L^Q$ are precisely the (smooth) solutions of (\ref{Aele}-\ref{Abdry1}-\ref{Abdry2}).

\item For every critical point $\gamma$ of $\mathbb{S}_L^Q$, the second Gateaux differential $d^2\mathbb{S}_L^Q(\gamma)$ of $\mathbb{S}_L^Q$ at $\gamma$ is Fredholm\footnote{ A continuous symmetric bilinear form $\alpha$ on the Hilbert space $\mathbb{H}$ is said to be Fredholm if the associated self-adjoint operator $A$ is Fredholm. Its Morse index is the dimension of the $A$-invariant subspace of $\mathbb{H}$ corresponding to the negative part of the spectrum of $\mathbb{H}$. So the Morse index of $\alpha$ is finite if and only if the negative spectrum of $A$ consists of finitely many eigenvalues with finite multiplicity. All of these notions do not depend on the choice of the Hilbert product on $\mathbb{H}$.} and has finite Morse index.

\end{enumerate}
\end{proposition}

\noindent\proof All the statements are of a local nature, so by
using the diffeomorphism $\varphi_*$ induced by an element
$(U,\varphi)$ of $\mathscr{F}$, we may assume that the smooth
Lagrangian is defined on $[0,1]\times U \times \R^n$, where $U$ is
an open subset of $\R^n$. By choosing $(U,\varphi)$ to be
bi-bounded, we may assume that $L$ satisfies (L1'), when proving
claim (i), and also (L2'), when proving claims (ii) and (iii).

If $\gamma\in H^1([0,1],U)$, $\xi\in H^1([0,1],\R^n)$, and
$\delta \in \R\setminus \{0\}$ has a small absolute value, we have
\begin{equation}
\label{regC1}
\begin{split}
\frac{1}{\delta} \Bigl(\mathbb{S}_L(\gamma+ \delta \xi) -\mathbb{S}_L
(\gamma)\Bigr) \\ = \int_0^1 \int_0^1 \Bigl( \frac{\partial L}{\partial q}
  (t,\gamma+ s \delta  \xi, \gamma'+ s \delta \xi') \cdot \xi
+ \frac{\partial L}{\partial v}(t,\gamma+ s \delta \xi,\gamma'+ s \delta \xi')
\cdot \xi' \Bigr) \,dt \, ds.
\end{split}
\end{equation}
The bounds (\ref{Agc1}) and the dominated convergence
theorem imply that the quantity (\ref{regC1}) converges to
\begin{equation}
\label{diff1}
D\mathbb{S}_L(\gamma)[\xi] := \int_0^1 \bigl( \frac{\partial L}{\partial q}
  (t,\gamma,\gamma')\cdot \xi + \frac{\partial L}{\partial v}(t,\gamma,\gamma')\cdot \xi' \bigr)\,
dt
\end{equation}
for $\delta \rightarrow 0$. Since $D\mathbb{S}_L(\gamma)$ is a bounded
linear functional on $H^1([0,1],\R^n)$, $\mathbb{S}_L$ is Gateaux
differentiable, and $D\mathbb{S}_L(\gamma)$ is its Gateaux
differential at $\gamma$.

In order to prove that $D\mathbb{S}_L$ is continuous at $\gamma$, we must show that $D\mathbb{S}_L(\gamma_h)$ converges to $D\mathbb{S}_L(\gamma)$ in the dual norm of $H^1$ when $\gamma_h$ converges to $\gamma$ in $H^1$. In particular, $\gamma_h$ converges to $\gamma$ uniformly, and there is a function $f$ in $L^2([0,1])$ such that for every $h\in \N$, $|\gamma_h'| \leq f$ almost everywhere (here we are using the fact that a sequence of real-valued functions which converges in $L^1$ is dominated almost-everywhere by an $L^1$ function).
By a standard argument involving subsequences\footnote{\label{subseq} That is, the fact that a sequence $(x_h)$ in a metric space converges to $x$ if and only if every subsequence of $(x_h)$ has a subsequence which converges to $x$.}, we may also assume that $\gamma_h' \rightarrow \gamma'$ almost everywhere.
Then the bounds (\ref{Agc1}) and the dominated convergence theorem imply that
\[
{\textstyle \frac{\partial L}{\partial q} (\cdot,\gamma_h,\gamma_h') \rightarrow \frac{\partial L}{\partial q} (\cdot,\gamma,\gamma') \mbox{ in } L^1([0,1]), \quad
\frac{\partial L}{\partial v} (\cdot,\gamma_h,\gamma_h') \rightarrow \frac{\partial L}{\partial v} (\cdot,\gamma,\gamma') \mbox{ in } L^2([0,1]).}
\]
These convergences imply that  $D\mathbb{S}_L(\gamma_h)$ converges to $D\mathbb{S}_L(\gamma)$ in the dual norm of $H^1$. Therefore,
$D\mathbb{S}_L(\gamma)$ depends continuously on $\gamma\in
H^1([0,1],U)$, so the total differential theorem implies that
$\mathbb{S}_L$ is continuously (Fr\'ech\'et) differentiable, and that
$D \mathbb{S}_L(\gamma)$ is its (Fr\'ech\'et) differential at $\gamma$.

Let $\gamma$ and $\xi$ be as above, and let $\eta$ be an element of $H^1([0,1],\R^n)$. By (L1') and by the dominated convergence theorem, the quantity
\begin{equation}
\label{regC2}
\begin{split}
\textstyle{
\frac{1}{\delta} \bigl(}& \textstyle{D\mathbb{S}_L(\gamma+ \delta \eta)[\xi] -
D\mathbb{S}_L(\gamma)[\xi]\bigr)} \\ =& \textstyle{
\int_0^1 \int_0^1 \Bigl( \frac{\partial^2 L}{\partial v^2}(t,
\gamma+ s \delta \eta, \gamma'+ s \delta \eta') \xi' \cdot \eta' +
\frac{\partial^2 L}{\partial v \partial q}
(t,\gamma+ s \delta \eta, \gamma'+ s \delta \eta')\xi' \cdot \eta} \\
& \textstyle{+ \frac{\partial^2 L}{\partial q \partial v}
(t,\gamma+ s\delta \eta, \gamma'+ s\delta \eta')\xi \cdot \eta'
+ \frac{\partial^2 L}{\partial q^2}
(t,\gamma+ s \delta\eta, \gamma'+ s \delta \eta')\xi \cdot \eta \Bigr) \, dt\,ds,}
\end{split}
\end{equation}
converges to
\begin{equation}
\label{diff2} \begin{split}
d^2 \mathbb{S}_L (\gamma)[\xi,\eta] :=
  \int_0^1 \Bigl( \frac{\partial^2 L}{\partial v^2} (t,\gamma,\gamma')
\xi' \cdot \eta' + \frac{\partial^2 L}{\partial v \partial q}
(t,\gamma,\gamma') \xi' \cdot \eta \\
+ \frac{\partial^2 L}{\partial q \partial v}
  (t,\gamma,\gamma')\xi \cdot \eta' + \frac{\partial^2 L}{\partial q^2}
(t,\gamma,\gamma') \xi \cdot \eta \Bigr) \, dt,
\end{split} \end{equation}
for $\delta \rightarrow 0$. Since $d^2 \mathbb{S}_L(\gamma)$ is a bounded
symmetric bilinear form on $H^1([0,1],\R^n)$, $D\mathbb{S}_L$ is Gateaux differentiable at $\gamma$, and its Gateaux differential at $\gamma$ is the bounded linear operator $D^2 \mathbb{S}_L(\gamma): H^1([0,1],\R^n) \rightarrow 
H^1([0,1],\R^n)^*$ defined by
\[
(D^2 \mathbb{S}_L(\gamma) \xi)[\eta] = d^2 \mathbb{S}_L(\gamma)[\xi,\eta], \quad \forall \xi,\eta \in H^1([0,1],\R^n).
\]
By (L1'), the map $\gamma \mapsto D^2 \mathbb{S}_L(\gamma)$ is bounded with respect to the norm-topology on the space of bounded self-adjoint operators, so the mean value theorem implies that $D\mathbb{S}_L$ is Lipschitz on convex subsets of $H^1([0,1],U)$. This concludes the proof of (i).

Let $Q$ be a submanifold of $M\times M$, and assume that $L$ satisfies also (L2). Let $\gamma$ be a critical point of $\mathbb{S}_L^Q$. By applying the above localization argument with $(U,\varphi)$ in $\mathscr{F}_Q$, we may assume that $\gamma$ is a critical point of $\mathbb{S}_L^W$, the restriction of $\mathbb{S}_L : H^1([0,1],U) \rightarrow \R$ to the intersection of the open set $H^1([0,1],U)$ with the closed linear subspace
\[
H^1_W([0,1],\R^n) = \set{ \gamma\in H^1([0,1],\R^n) }{(\gamma(0),\gamma(1))\in W},
\]
where $W$ is the linear subspace of $\R^n \times \R^n$ such that $(\varphi(0,q),\varphi(1,q))\in Q$ for every $q\in U\cap W$. The boundary condition (\ref{Abdry1}) is equivalent to
\[
(\gamma(0),\gamma(1)) \in W,
\]
and it is satisfied by every element of $H^1_W([0,1],\R^n)$. Differentiating the condition $(\varphi(0,q),\varphi(1,q))\in Q$ for every $q\in U\cap W$ at the point $(q_0,q_1)\in U \cap W$, we obtain that the linear map
\[
(\xi_0,\xi_1) \mapsto (D_q \varphi(0,q_0)[\xi_0], D_q \varphi(1,q_1)[\xi_1])
\]
maps $W$ isomorphically onto the tangent space of $Q$ at the point $(\varphi(0,q_0),\varphi(1,q_1))$, so the boundary condition (\ref{Abdry2}) is equivalent to
\begin{equation}
\label{bdry2}
\frac{\partial L}{\partial v}(0,\gamma(0),\gamma'(0)) \cdot \xi_0 = \frac{\partial L}{\partial v}(1,\gamma(1),\gamma'(1)) \cdot \xi_1, \quad \forall (\xi_0,\xi_1) \in W.
\end{equation}
For every smooth curve $\xi:[0,1]\rightarrow \R^n$ which has compact support in $]0,1[$, formula (\ref{diff1}) and an integration by parts produce the identity
\[
\int_0^1 \left( - \int_0^t \frac{\partial L}{\partial q}(s,\gamma(s),\gamma'(s))\, ds + \frac{\partial L}{\partial v} (t,\gamma(t),\gamma'(t)) \right) \cdot \xi'(t) \, dt = 0.
\]
Then Du Bois-Reymond Lemma implies that there exists a vector $u\in \R^n$ such that
\begin{equation}
\label{impl}
- \int_0^t \frac{\partial L}{\partial q}(s,\gamma(s),\gamma'(s))\, ds + \frac{\partial L}{\partial v} (t,\gamma(t),\gamma'(t)) = u \quad \mbox{a.e. in } [0,1].
\end{equation}
By (L2'), the map
\begin{equation}
\label{diffeo}
[0,1] \times U \times \R^n \rightarrow [0,1] \times U \times \R^n, \quad (t,q,v) \mapsto \left(t,q,\frac{\partial L}{\partial v} (t,q,v) \right),
\end{equation}
is a surjective smooth diffeomorphism. We denote by
\[
(t,q,p) \rightarrow (t,q,\psi(t,q,p))
\]
its inverse. Then (\ref{impl}) is equivalent to
\[
\gamma'(t) = \psi\left(t,\gamma(t), u + \int_0^t \frac{\partial L}{\partial q}(s,\gamma(s),\gamma'(s))\, ds \right)\quad \mbox{a.e. in } [0,1],
\]
and a boot-strap argument shows that $\gamma$ is smooth. We can then apply a different integration by parts to the identity $D\mathbb{S}_L^W(\gamma)[\xi]=0$, and we obtain that the quantity
\begin{equation}
\label{parts}
\begin{split}
\int_0^1 \left( \frac{\partial L}{\partial q}(t,\gamma,\gamma') - \frac{d}{dt} \Bigl( \frac{\partial L}{\partial v} (t,\gamma,\gamma') \Bigr) \right)\, dt \\ + \frac{\partial L}{\partial v}(1,\gamma(1),\gamma'(1)) \cdot \xi(1) - \frac{\partial L}{\partial v}(0,\gamma(0),\gamma'(0)) \cdot \xi(0)
\end{split}
\end{equation}
vanishes
for every smooth curve $\xi:[0,1]\rightarrow \R^n$ such that the pair $(\xi(0),\xi(1))$ belongs to $W$.
By taking curves $\xi$ with compact support we obtain that $\gamma$ must solve (\ref{Aele}). Then, letting $\xi$ vary among all smooth curves $\xi$ such that $(\xi(0),\xi(1))\in W$, we find that (\ref{bdry2}) holds.

 This shows that every critical point of $\mathbb{S}_L^Q$ is a smooth solution of (\ref{Aele}-\ref{Abdry1}-\ref{Abdry2}). Conversely, the fact that (\ref{diffeo}) is a diffeomorphism and the differentiable dependence of solutions of ordinary differential equations on the coefficients imply that every solution of (\ref{Aele}) is smooth. If the boundary conditions (\ref{Abdry1}) and (\ref{Abdry2}) are also satisfied, by integrating by parts the expression for $D\mathbb{S}_L^Q(\gamma)[\xi]$ as in (\ref{parts}), one immediately sees that $\gamma$ is a critical point of $\mathbb{S}_L^Q$.
This concludes the proof of (ii).

Let $\gamma$ be a critical point of $\mathbb{S}_L^Q$. Then the Gateaux second differential of $\mathbb{S}_L^Q$ at $\gamma$ is well-defined as a symmetric continuous bilinear form on $T_{\gamma} H^1_Q([0,1],M)$. By using the above localization argument, we may identify $\gamma$ with a critical point of $\mathbb{S}_L^W$ in $H^1_W([0,1],\R^n)$, and we may identify $d^2 \mathbb{S}_L^Q(\gamma)$ with $d^2 \mathbb{S}_L^W(\gamma)$, that is the restriction of the symmetric bilinear form (\ref{diff2}) to the Hilbert space $H^1_W([0,1],\R^n)$.
By (L2'), the self-adjoint operator on $H^1_W([0,1],\R^n)$ representing the symmetric bilinear form
\[
\alpha(\gamma)[\xi,\eta] := \int_0^1 \frac{\partial^2 L}{\partial v^2} (t,\gamma,\gamma') \xi' \cdot \eta' \, dt
\]
with respect to the Hilbert product is Fredholm and non-negative. The remaining three terms in (\ref{diff2}) are bilinear forms which are continuous on $H^1\times L^2$, $L^2 \times H^1$, and $L^2 \times L^2$, respectively. The compactness of the embedding $H^1\hookrightarrow L^2$ then implies that the self-adjoint operator representing $d^2 \mathbb{S}_L^Q(\gamma) - \alpha(\gamma)$ is compact. This implies (iii).
\qed

The Lagrangian action functional $\mathbb{S}_L$ is not of class $C^2$, unless $L$ is a polynomial of degree at most two on each fiber of $TM$. In fact, the first term in (\ref{diff2}), that is, the symmetric bilinear form $\alpha(\gamma)$ defined above, depends continuously on $\gamma\in H^1$ if and only if $\gamma \mapsto \partial^2 L/\partial v^2 (t,\gamma,\gamma')$ is continuous from $H^1$ to $L^{\infty}$, and the latter fact is true if and only if $\partial^2 L/\partial v^2 (t,q,v)$ does not depend on $v$ (for $(t,q)$ in a neighborhood of $(t,\gamma(t))$). This is a manifestation of a well-known phenomenon concerning the differentiability of Nemitsky operators on $L^2$, explained for instance in \cite{ap95} (see in particular Proposition 2.8 of Chapter 1). Actually, in general the functional $\mathbb{S}_L$ also fails to be twice differentiable, as shown by the following:

\begin{proposition}
\label{contro}
Assume that the smooth Lagrangian function $L:[0,1]\times TM \rightarrow \R$ satisfies (L1). If $\mathbb{S}_L:H^1([0,1],M)\rightarrow \R$ is twice differentiable at the curve $\gamma\in H^1([0,1],M)$, then for every $t\in [0,1]$ the function
\[
T_{\gamma(t)} M \rightarrow \R, \quad v\mapsto L(t,\gamma(t),v),
\]
is a polynomial of degree at most two.
\end{proposition}

\noindent\proof By using a suitable coordinate system
$(U,\varphi)\in \mathscr{F}$, we may assume that $L$ is defined on
$[0,1]\times U \times \R^n$, where $U$ is an open subset of $\R^n$,
that it satisfies (L1'), and that $\gamma\in H^1([0,1],U)$. Since
$\gamma$ is continuous and $\gamma'$ exists a.e., the thesis is
equivalent to the fact that for almost every $t\in [0,1]$ there
holds
\begin{equation}
\label{tesiA}
\frac{\partial L}{\partial v} (t,\gamma(t),\gamma'(t)+v) - \frac{\partial L}{\partial v} (t,\gamma(t),\gamma'(t)) - \frac{\partial^2 L}{\partial v^2} (t,\gamma(t),\gamma'(t)) v = 0,
\end{equation}
for every $v\in \R^n$.
Since $\partial^2 L/\partial v^2$ is bounded, the map appearing in the left-hand side of the above identity is measurable with respect to $t$ and Lipschitz in $v$, uniformly with respect to $t$. Under these assumptions, it is easy to show that the fact that for almost every $t\in [0,1]$ (\ref{tesiA}) holds is equivalent to the fact that for every $v\in \R^n$ the quantity appearing in the left-hand side of (\ref{tesiA}) vanishes for almost every $t\in [0,1]$
(actually, this fact is true in the more general case of Caratheodory
functions).
So, striving for a contradiction, we can assume that there is a set of positive measure $J\subset [0,1]$, two non-zero vectors $v,w\in \R^n$, and a positive number $c$ such that
\begin{equation}
\label{contra}
\left( \frac{\partial L}{\partial v} (t,\gamma(t),\gamma'(t)+v) - \frac{\partial L}{\partial v} (t,\gamma(t),\gamma'(t)) - \frac{\partial^2 L}{\partial v^2} (t,\gamma(t),\gamma'(t)) v \right) \cdot w > c,
\end{equation}
for every $t\in J$.
For every $\epsilon>0$ smaller than the measure of $J$ we choose a subset $J_{\epsilon}\subset J$ of measure $\epsilon$, in such a way that $J_{\epsilon} \subset J_{\epsilon'}$ if $\epsilon<\epsilon'$. We define elements $\eta_{\epsilon}$ and $\xi_{\epsilon}$ of $H^1([0,1],\R^n)$ by
\[
\eta_{\epsilon} (t) = v \int_0^t I_{J_{\epsilon}}(s)\, ds, \quad
\xi_{\epsilon} (t) = w \int_0^t I_{J_{\epsilon}}(s)\, ds,
\]
where $I_{J_{\epsilon}}$ is the characteristic function of $J_{\epsilon}$.
Then, using the Hilbert norm $\|u\|_{H^1} = (|u(0)|^2 + \|u'\|_{L^2}^2)^{1/2}$
on $H^1([0,1],\R^n)$, we have
\[
\|\eta_{\epsilon} \|_{H^1} = \left( \int_{J_{\epsilon}} |v|^2\, dt \right)^{1/2}=
|v| \epsilon^{1/2}, \quad \|\xi_{\epsilon} \|_{H^1} = |w| \epsilon^{1/2}.
\]
By assumption, $D\mathbb{S}_L$ is differentiable at $\gamma$ and, by
Proposition \ref{reg} (i), its differential at $\gamma$ coincides with the Gateaux differential $D^2 \mathbb{S}_L(\gamma)$. In particular,
\begin{equation}
\label{asint}
D\mathbb{S}_L(\gamma + \eta_{\epsilon})[\xi_{\epsilon}] - D\mathbb{S}_L(\gamma)[\xi_{\epsilon}] - d^2 \mathbb{S}_L(\gamma)[\xi_{\epsilon},\eta_{\epsilon}] = o(\|\eta_{\epsilon}\|_{H^1}) \|\xi_{\epsilon}\|_{H^1} = o(\epsilon),
\end{equation}
for $\epsilon \rightarrow 0$.
If we express the left-hand side of (\ref{asint}) by formulas (\ref{diff1}) and (\ref{diff2}), we deduce that
\begin{equation}
\label{mostro} \begin{split}
\int_0^1 \Bigl( \frac{\partial L}{\partial q} (t,\gamma+\eta_{\epsilon},\gamma'+\eta_{\epsilon}') - \frac{\partial L}{\partial q} (t,\gamma,\gamma') - \frac{\partial^2 L}{\partial q^2} (t,\gamma,\gamma') \eta_{\epsilon} \\ - \frac{\partial^2 L}{\partial v \partial q} (t,\gamma,\gamma') \eta_{\epsilon}' \Bigr)\cdot \xi_{\epsilon} \, dt + \int_0^1 \Bigl( \frac{\partial L}{\partial v} (t,\gamma+\eta_{\epsilon},\gamma'+\eta_{\epsilon}') - \frac{\partial L}{\partial v} (t,\gamma,\gamma')  \\ - \frac{\partial^2 L}{\partial q \partial v} (t,\gamma,\gamma') \eta_{\epsilon} - \frac{\partial^2 L}{\partial v^2} (t,\gamma,\gamma') \eta_{\epsilon}' \Bigr)\cdot \xi_{\epsilon}' \, dt  = o(\epsilon),
\end{split} \end{equation}
for $\epsilon \rightarrow 0$.
Since $\|\xi_{\epsilon}\|_{\infty} = \epsilon |w|$, the absolute value of the first integral in (\ref{mostro}) can be estimated from above by
\[
\epsilon |w| \int_0^1 \left| \frac{\partial L}{\partial q} (t,\gamma+\eta_{\epsilon},\gamma'+\eta_{\epsilon}') - \frac{\partial L}{\partial q} (t,\gamma,\gamma) - \frac{\partial^2 L}{\partial q^2} (t,\gamma,\gamma') \eta_{\epsilon} - \frac{\partial^2 L}{\partial v \partial q} (t,\gamma,\gamma') \eta_{\epsilon}' \right| \, dt.
\]
Since $\eta_{\epsilon}$ and $\eta_{\epsilon}'$ converge to zero almost everywhere, by (L1') the dominated convergence theorem implies that the above integral tends to zero for $\epsilon\rightarrow 0$. Therefore,
\begin{equation}
\label{primo}
\begin{split}
\int_0^1 \Bigl( \frac{\partial L}{\partial q} (t,\gamma+\eta_{\epsilon},\gamma'+\eta_{\epsilon}') - \frac{\partial L}{\partial q} (t,\gamma,\gamma') - \frac{\partial^2 L}{\partial q^2} (t,\gamma,\gamma') \eta_{\epsilon} \\ - \frac{\partial^2 L}{\partial v \partial q} (t,\gamma,\gamma') \eta_{\epsilon}' \Bigr)\cdot \xi_{\epsilon} \, dt  = o(\epsilon),
\end{split}
\end{equation}
for $\epsilon\rightarrow 0$. By the definition of $\eta_{\epsilon}$ and $\xi_{\epsilon}$, we have
\begin{eqnarray*}
\int_0^1 \left( \frac{\partial L}{\partial v} (t,\gamma+\eta_{\epsilon},\gamma'+\eta_{\epsilon}') - \frac{\partial L}{\partial v} (t,\gamma,\gamma'+\eta_{\epsilon}') - \frac{\partial^2 L}{\partial q \partial v} (t,\gamma,\gamma'+ \eta_{\epsilon}') \eta_{\epsilon} \right)\cdot \xi_{\epsilon}' \, dt \\
= \int_{J_{\epsilon}} \left( \frac{\partial L}{\partial v} (t,\gamma+\eta_{\epsilon},\gamma'+v) - \frac{\partial L}{\partial v} (t,\gamma,\gamma'+v) - \frac{\partial^2 L}{\partial q \partial v} (t,\gamma,\gamma'+v) \eta_{\epsilon} \right)\cdot w \, dt \\
= \int_{J_{\epsilon}} \int_0^1 \left( \frac{\partial^2 L}{\partial q \partial v} (t,\gamma+s\eta_{\epsilon},\gamma'+v) - \frac{\partial^2 L}{\partial q \partial v} (t,\gamma,\gamma'+v) \right) \eta_{\epsilon} \cdot w \, ds \, dt.
\end{eqnarray*}
Since $\|\eta_{\epsilon}\|_{\infty} = \epsilon |v|$, the absolute value of the last double integral can be estimated from above by
\begin{eqnarray*}
\epsilon |v| |w| \int_{J_{\epsilon}} \int_0^1 \left| \frac{\partial^2 L}{\partial q \partial v} (t,\gamma+s\eta_{\epsilon},\gamma'+v) - \frac{\partial^2 L}{\partial q \partial v} (t,\gamma,\gamma'+v) \right| \, ds \, dt \\
\leq \epsilon |v| |w| \epsilon^{1/2} \left(
\int_0^1 \int_0^1 \Bigl| \frac{\partial^2 L}{\partial q \partial v} (t,\gamma+s\eta_{\epsilon},\gamma'+v) - \frac{\partial^2 L}{\partial q \partial v} (t,\gamma,\gamma'+v) \Bigr|^2 \, ds \, dt \right)^{1/2},
\end{eqnarray*}
where we have also used the Cauchy-Schwarz inequality. By (L1'), the above double integral is bounded, so we have shown that
\begin{equation}
\label{seconda}
\begin{split}
\int_0^1 \Bigl( \frac{\partial L}{\partial v} (t,\gamma+\eta_{\epsilon},\gamma'+\eta_{\epsilon}') - \frac{\partial L}{\partial v} (t,\gamma,\gamma'+\eta_{\epsilon}') \\ - \frac{\partial^2 L}{\partial q \partial v} (t,\gamma,\gamma'+ \eta_{\epsilon}') \eta_{\epsilon} \Bigr)\cdot \xi_{\epsilon}' \, dt = O(\epsilon^{3/2}),
\end{split} \end{equation}
for $\epsilon\rightarrow 0$. Similarly,
\begin{eqnarray*}
\left| \int_0^1 \Bigl( \frac{\partial^2 L}{\partial q \partial v} (t,\gamma,\gamma'+\eta_{\epsilon}') -\frac{\partial^2 L}{\partial q \partial v} (t,\gamma,\gamma') \Bigr) \eta_{\epsilon} \cdot \xi_{\epsilon}'\, dt \right| \\
= \left| \int_{J_{\epsilon}} \Bigl( \frac{\partial^2 L}{\partial q \partial v} (t,\gamma,\gamma'+v) -\frac{\partial^2 L}{\partial q \partial v} (t,\gamma,\gamma') \Bigr) \eta_{\epsilon} \cdot w\, dt \right| \\
\leq \epsilon |v| |w| \int_{J_{\epsilon}} \left| \frac{\partial^2 L}{\partial q \partial v} (t,\gamma,\gamma'+v) -\frac{\partial^2 L}{\partial q \partial v} (t,\gamma,\gamma') \right| \, dt \\
\leq \epsilon |v| |w| \epsilon^{1/2} \left( \int_0^1 \Bigl| \frac{\partial^2 L}{\partial q \partial v} (t,\gamma,\gamma'+v) -\frac{\partial^2 L}{\partial q \partial v} (t,\gamma,\gamma') \Bigr|^2 \, dt \right)^{1/2}.
\end{eqnarray*}
Since the last integral is bounded because of (L1'), we have
\begin{equation}
\label{terza}
\int_0^1 \Bigl( \frac{\partial^2 L}{\partial q \partial v}
(t,\gamma,\gamma'+\eta_{\epsilon}') -\frac{\partial^2 L}{\partial q
\partial v} (t,\gamma,\gamma') \Bigr) \eta_{\epsilon}
\cdot \xi_{\epsilon}'\, dt = O(\epsilon^{3/2}),
\end{equation}
for $\epsilon \rightarrow 0$. By comparing (\ref{mostro}) with (\ref{primo}), (\ref{seconda}), and (\ref{terza}), we deduce that
\begin{equation}
\label{resto}
\int_0^1 \left( \frac{\partial L}{\partial v} (t,\gamma,\gamma'+\eta_{\epsilon}') - \frac{\partial L}{\partial v} (t,\gamma,\gamma') - \frac{\partial^2 L}{\partial v^2} (t,\gamma,\gamma') \eta_{\epsilon}' \right)\cdot \xi_{\epsilon}' \, dt = o(\epsilon), \quad \mbox{for } \epsilon \rightarrow 0.
\end{equation}
However, the left-hand side of (\ref{resto}) equals
\[
\int_{J_{\epsilon}} \left( \frac{\partial L}{\partial v} (t,\gamma,\gamma'+v) - \frac{\partial L}{\partial v} (t,\gamma,\gamma') - \frac{\partial^2 L}{\partial v^2} (t,\gamma,\gamma') v \right)\cdot w \, dt,
\]
which by (\ref{contra}) is larger than $c \epsilon$. This contradiction concludes the proof.
\qed

In particular, a Lagrangian function $L$ induces a twice differentiable action functional on $H^1([0,1],M)$ if and only if it is of the form (\ref{em}). In this case, it is easy to show that the action functional is actually smooth on $H^1([0,1],M)$. In particular, electro-magnetic Lagrangian functions - that is functions of the form (\ref{em}) with $A$ positive - are the only Lagrangian functions which satisfy (L2) and induce a smooth action functional.

\begin{remark}
{\rm Similarly, if the energy functional induced by a Finsler metric
on $M$ is twice differentiable at a curve $\gamma\in H^1([0,1],M)$
whose velocity is different from zero a.e., then the Finsler metric
is actually Riemannian along $\gamma$. This fact can be shown by
adapting the argument of Proposition \ref{contro}, taking into
account the fact in this case the Lagrangian function, that is the
Finsler energy, is not twice differentiable at the zero section of
$TM$. See \cite{rad92} and \cite{cjm09} for two different ways of
developing a Morse theory in the Finsler setting.}
\end{remark}

\paragraph{The Palais-Smale condition.} We conclude this section by proving the Palais-Smale condition for the functional $\mathbb{S}_L^Q$.
The metric $g$ on $M$ induces a Riemannian metric on the Hilbert manifold $H^1([0,1],M)$, namely
\[
\langle \xi, \eta \rangle := \int_0^1 g(\nabla_t \xi, \nabla_t \eta) \, dt + \int_0^1 g(\xi,\eta) \, dt, \quad \begin{array}{l} \forall \gamma \in H^1([0,1],M), \\ \forall \xi,\eta \in T_{\gamma} H^1([0,1],M), \end{array}
\]
where $\nabla_t$ denotes the Levi-Civita covariant derivative along $\gamma$ induced by $g$. The corresponding norm on the tangent bundle of $H^1([0,1],M)$ is denoted by $\|\cdot\|$. If $(U,\varphi)\in \mathscr{F}$ is bi-bounded,
then the Riemannian metric $\langle \cdot , \cdot \rangle$ on $\mathscr{U}_{(U,\varphi)}$ is uniformly equivalent to the one induced from the Hilbert structure of $H^1([0,1],\R^n)$ by $\varphi_*$.

The Riemannian metric $\langle \cdot, \cdot \rangle$ is complete if $g$ is complete. This can be shown by proving that every Cauchy sequence $(\gamma_h)$ for the distance function induced by $\langle \cdot, \cdot \rangle$ is a fortiori a Cauchy sequence for the uniform distance on $C^0([0,1],M)$ induced by $g$. Then, since $C^0([0,1],M)$ is complete and the sets $\mathscr{U}_{(U,\varphi)}$ are $C^0$-open, we may assume that $(\gamma_h)$ is a Cauchy sequence in $H^1([0,1],U)$, with $U$ an open subset of $\R^n$, and that it converges uniformly to some continuous curve $\gamma:[0,1]\rightarrow U$. Then the completeness of $H^1([0,1],\R^n)$ implies that $\gamma_h \rightarrow \gamma$ in $H^1$.
Therefore $(H^1([0,1],M), \langle \cdot,\cdot \rangle)$ is complete.

Since $H^1_Q([0,1],M)$ is closed in $H^1([0,1],M)$ when $Q$ is closed in
$M\times M$, the latter fact implies that the  restriction of the metric $\langle \cdot,\cdot \rangle$ to the submanifold $H^1_Q([0,1],M)$ is also complete.
The gradient vector field of $\mathbb{S}_L^Q$ with respect to this Riemannian metric is denoted by $\grad \mathbb{S}_L^Q$.

The following result is due to Benci, \cite{ben86}. The proof we present here is taken from Appendix A in \cite{af07}.

\begin{proposition}
\label{ps}
Assume that $(M,g)$ is complete, that $Q$ is closed in $M\times M$, and that at least one of the sets $\pi_1(Q)$ and $\pi_2(Q)$ is bounded, where $\pi_1$ and $\pi_2$ denote the two projections $M\times M \rightarrow M$. Assume that $L$ satisfies (L1), (L2), and
\begin{equation}
\label{unif}
L(t,q,v) \geq \ell_0 \, g(v,v) - c, \quad \forall (t,q,v) \in [0,1]\times TM,
\end{equation}
where $\ell_0 >0$ and $c\in \R$. Then the pair $(- \grad \mathbb{S}_L^Q, \mathbb{S}_L^Q)$ satisfies the Palais-Smale condition on $H^1_Q([0,1],M)$. In other words, the functional $\mathbb{S}_L^Q$
satisfies the Palais-Smale condition
with respect to the metric $\langle \cdot,\cdot \rangle$ in the usual sense.
\end{proposition}

\noindent\proof Let $(\gamma_h)$ be a sequence in $H^1_Q([0,1],M)$
such that $\mathbb{S}_L^Q(\gamma_h)$ is bounded and $\| \grad
\mathbb{S}_L^Q(\gamma_h)\|$ is infinitesimal. The estimate
(\ref{unif}) and the upper bound on $\mathbb{S}_L^Q(\gamma_h)$
\newpage \noindent imply that the sequence
\begin{equation}
\label{bdd}
\int_0^1 g(\gamma_h',\gamma_h')\, dt
\end{equation}
is bounded. By the Cauchy-Schwarz inequality, if $0\leq s \leq t \leq 1$ we have
\[
\dist (\gamma_h(t), \gamma_h(s)) \leq \int_s^t g(\gamma_h',\gamma_h')^{1/2} \, d\sigma \leq |s-t|^{1/2} \left( \int_0^1 g(\gamma_h',\gamma_h')\, d\sigma \right)^{1/2},
\]
so $(\gamma_h)$ is equi-1/2-H\"older continuous. Moreover, by the assumption on $Q$ at least one of the sequences $(\gamma_h(0))$ and $(\gamma_h(1))$ is bounded. Since $(M,g)$ is complete, the Ascoli-Arzel\`a theorem implies that the sequence $(\gamma_h)$ is compact in the $C^0$ topology, so, up to the choice of a subsequence,  we may assume that $(\gamma_h)$ converges to some continuous curve $\gamma:[0,1]\rightarrow M$, with $(\gamma(0),\gamma(1))\in Q$, because $Q$ is closed. Since
\[
\set{\mathscr{U}_{(U,\varphi)}}{(U,\varphi)\in \mathscr{F}_Q \mbox{ bi-bounded}}
\]
is a covering of $H^1_Q([0,1],M)$ consisting of $C^0$-open sets, we may also assume that $\gamma_h$
belongs to the same $\mathscr{U}_{(U,\varphi)}$ for every $h\in \N$, for some $(U,\varphi)\in \mathscr{F}_Q$ bi-bounded and such that $\gamma$ belongs to the image of $\varphi_*$ extended to the space of continuous $U$-valued curves.

By using the induced diffeomorphism $\varphi_*$, we can localize the analysis to $\R^n$ and we may assume that the Lagrangian function $L$ is defined on $[0,1]\times U \times \R^n$, with $U$ an open neighborhood of $0$, and that $L$ satisfies (L1') and (L2'). There is a linear subspace $W$ of $\R^n\times \R^n$ such that the sequence $(\gamma_h)$ belongs to
\[
H^1_W([0,1],U) = \set{\gamma\in H^1([0,1],\R^n)}{\gamma([0,1])\subset U, \; (\gamma(0),\gamma(1)) \in W}.
\]
Since the Riemannian metric $\langle \cdot,\cdot \rangle$ on $\mathscr{U}_{(U,\varphi)}$ is uniformly equivalent to the one induced from the Hilbert product of $H^1_W([0,1],\R^n)$ by $\varphi_*$, the sequence $(\gamma_h)$ satisfies
\[
D\mathbb{S}_L^W (\gamma_h) \rightarrow 0,
\]
in the norm topology of the dual of the Hilbert space $H^1_W([0,1],\R^n)$.
Moreover, $(\gamma_h)$
converges uniformly to some continuous curve $\gamma:[0,1]\rightarrow U$ such that $(\gamma(0),\gamma(1))\in W$, and it is bounded in $H^1_W([0,1],\R^n)$ (because (\ref{bdd}) is bounded). Under these assumptions, we must prove that $(\gamma_h)$ has a subsequence which converges in the Hilbert space topology to an element of $H^1_W([0,1],U)$.

Since $(\gamma_h)$ is bounded in $H^1_W([0,1],\R^n)$, up to the choice of a subsequence we may assume that it converges weakly in $H^1$ and uniformly to some element $\tilde{\gamma}$ of $H^1_W([0,1],\R^n)$. Since $(\gamma_h)$ converges uniformly to $\gamma$, we have that $\tilde{\gamma}= \gamma$, so $\gamma$ belongs to $H^1_W([0,1],U)$. It remains to show that the convergence of $(\gamma_h)$ to $\gamma$ is strong in $H^1$.

The sequence $D\mathbb S_L^W(\gamma_h)[\gamma_h-\gamma]$ is
infinitesimal, that is
\[
\int_0^1  \frac{\partial L}{\partial q}(t,\gamma_h,\gamma_h') \cdot
(\gamma_h-\gamma)\, dt  + \int_0^1 \frac{\partial
L}{\partial v}(t,\gamma_h,\gamma_h')\cdot (\gamma_h'-\gamma')
\, dt \rightarrow 0 \quad \mbox{for } h \rightarrow \infty.
\]
Since $\partial L/\partial q(\cdot ,\gamma_h,\gamma_h')$ is bounded
in $L^1$ by (\ref{Agc1}), and $\gamma_h-\gamma$ converges to $0$
in $L^{\infty}$, the first integral in the above expression tends to zero, so
\begin{equation}
\label{inf}
\int_0^1 \frac{\partial L}{\partial v}(t,\gamma_h,\gamma_h') \cdot (\gamma_h'-\gamma')
\, dt \rightarrow 0 \quad \mbox{for } h \rightarrow \infty.
\end{equation}
For a.e.\ $t\in [0,1]$ we have by (L2')
\begin{eqnarray*}
\frac{\partial L}{\partial v} (t,\gamma_h,\gamma_h')\cdot (\gamma_h' - \gamma') - \frac{\partial L}{\partial v} (t,\gamma_h,\gamma') \cdot (\gamma_h' - \gamma') \\ = \int_0^1 \frac{\partial^2 L}{\partial v^2}
 (t,\gamma_h,\gamma' + s(\gamma_h' - \gamma')) (\gamma_h' -
\gamma') \cdot (\gamma_h' - \gamma') \, ds   \geq \ell_2 |\gamma_h'(t) -
\gamma'(t)|^2.
\end{eqnarray*}
Integrating this inequality over $[0,1]$ we get
\begin{equation}
\label{dis}
\begin{split}
\ell_2 \int_0^1 |\gamma_h' - \gamma'|^2 \, dt \\ \leq \int_0^1 \frac{\partial L}{\partial v} (t,\gamma_h,\gamma_h') \cdot (\gamma_h'-\gamma') \, dt -
\int_0^1  \frac{\partial L}{\partial v} (t,\gamma_h,\gamma') \cdot (\gamma_h'-\gamma') \, dt.
\end{split}
\end{equation}
The first integral on the right-hand side tends to $0$ for
$h\rightarrow \infty$ by (\ref{inf}). Since $\partial L/\partial v
(\cdot,\gamma_h,\gamma')$ converges strongly to $\partial L/\partial v
(\cdot,\gamma,\gamma')$ in $L^2$ by (\ref{Agc1}), and since
$(\gamma_h'-\gamma')$ converges weakly to $0$ in $L^2$, the last
integral in (\ref{dis}) is also infinitesimal. Therefore,
(\ref{dis}) implies that $(\gamma_h)$ converges strongly to
$\gamma$ in $H^1$, concluding the proof.
\qed
\begin{remark}
{\rm The completeness of the metric $g$ is associated to the uniform
lower estimate (\ref{unif}). In fact, as already observed,
conditions (L1) and (L2) are semi-local and do not depend on the
choice of the metric. Notice also that if $M$ is compact,
(\ref{unif}) is implied by (L2).}
\end{remark}
\begin{remark}
\label{bdgrad} {\rm Under the assumptions of Proposition \ref{ps},
the gradient of $\mathbb{S}_L^Q$ is bounded on every sublevel of
$\mathbb{S}_L^Q$. In fact, the argument used at the beginning of the
proof shows that every sublevel of $\mathbb{S}_L^Q$ is pre-compact
in $C^0$. As such, it can be covered by finitely many sets
$\mathscr{U}_{(\varphi,U)}$, with $(U,\varphi)\in \mathscr{F}_Q$
bi-bounded. Then the boundedness of $\|\grad \mathbb{S}_L^Q\|$ can
be checked locally, by showing that $D\mathbb{S}_L^W$ is bounded on
bounded subsets of $H^1_W([0,1],U)$. The latter fact follows from
the formula (\ref{diff1}) together with the bounds (\ref{Agc1}).}
\end{remark}

\section{\bf A smooth pseudo-gradient}
\def\theequation{4.\arabic{equation}}\makeatother
\setcounter{equation}{0}

Let us assume that all the solutions of  the problem (\ref{Aele}-\ref{Abdry1}-\ref{Abdry2}) are non-degenerate, meaning that:

\begin{enumerate}
\item[(L0)] The linear boundary value problem obtained by linearizing (\ref{Aele}-\ref{Abdry1}-\ref{Abdry2}) along any solution has only the trivial solution.
\end{enumerate}

Since the solutions of the  linear boundary value problem obtained by linearizing (\ref{Aele}-\ref{Abdry1}-\ref{Abdry2}) along the solution $\gamma$ are precisely the elements of the kernel of $d^2 \mathbb{S}_L^Q(\gamma)$,
assumption (L0) is equivalent to the fact that every critical point of $\mathbb{S}_L^Q$ is non-degenerate.
The  aim of this section is to prove the following:

\begin{theorem}
\label{lyap}
Assume that $(M,g)$ is complete, that $Q$ is closed in $M\times M$, and that at least one of the sets $\pi_1(Q)$ and $\pi_2(Q)$ is bounded, where $\pi_1$ and $\pi_2$ denote the two projections $M\times M \rightarrow M$. Assume that $L$ satisfies (L0) with respect to $Q$, (L1), (L2), and
\[
L(t,q,v) \geq \ell_0 \, g(v,v) - c, \quad \forall (t,q,v) \in [0,1]\times TM,
\]
where $\ell_0 >0$ and $c\in \R$.
Then, there exists a smooth vector field $X$ on the Hilbert manifold $H^1_Q([0,1],M)$ such that:
\begin{enumerate}
\item $X$ is complete.
\item $\mathbb{S}_L^Q$ is a Lyapunov function for $X$.
\item $X$ is a Morse vector field. More precisely,  $\gamma$ is a critical point of $\mathbb{S}_L^Q$ if and only if $X(\gamma)=0$, and $\nabla X(\gamma) = - D^2 \mathbb{S}_L^Q(\gamma)$, where $D^2 \mathbb{S}_L^Q(\gamma)$ is the self-adjoint automorphism of the tangent space of $H^1_Q([0,1],M)$ at $\gamma$ which represents the second Gateaux differential of $\mathbb{S}_L^Q$ at $\gamma$ with respect to some equivalent Hilbert product.
\item The pair $(X,\mathbb{S}_L^Q)$ satisfies the Palais-Smale condition.
\item $X$ satisfies the Morse-Smale condition up to every order.
\end{enumerate}
\end{theorem}

The main point in the proof of this result is to construct the vector field $X$ in a neighborhood of each critical point of $\mathbb{S}_L^Q$.
This construction uses some facts about the strong convergence in the algebra of bounded operators on the real Hilbert space $\mathbb{H}$, that we now recall. A sequence of bounded operators $(A_h)$ is said to converge {\em strongly} to the bounded operator $A$ if $A_h \xi$ converges to $A \xi$ for every vector $\xi$ in $\mathbb{H}$. By the principle of uniform boundedness, strongly convergent sequences of operators are bounded in norm.
If $(A_h)$ converges strongly to $A$, and $(K_h)$ is a sequence of compact operators which converges in norm to the (necessarily compact) operator $K$, then it is easy to show that $(A_h K_h)$ converges to $AK$ in norm. Furthermore, the functional calculus on the space of self-adjoint operators on $\mathbb{H}$ is strongly continuous, meaning that if the sequence of self-adjoint bounded operators $(A_h)$ converges to $A$ strongly and $f$ is a continuous real function on a closed subset of $\R$ containing the spectrum of every $A_h$ and of $A$, then $f(A_h)$ converges to $f(A)$ strongly. This fact can be proved by
approximating $f$ with polynomials, see e.g.\ Proposition 5.3.2
in \cite{kr97}.

\begin{lemma}
\label{loc}
Let $L: [0,1] \times \R^n \times \R^n \rightarrow \R$ be a smooth Lagrangian function satisfying (L1') and (L2'). We also assume that $0$ is a non-degenerate critical point of $\mathbb{S}_L^W$, the restriction of $\mathbb{S}_L$ to the Hilbert space $H^1_W([0,1],\R^n)$, for a given linear subspace $W$ of $\R^n \times \R^n$. Then there exists an equivalent Hilbert product $\langle \cdot , \cdot \rangle_0$ on $H^1_W([0,1],\R^n)$, with associated norm $\|\cdot\|_0$, and a positive number $\lambda$ such that, denoting by $D^2 \mathbb{S}^W_L(0)$ the self-adjoint operator representing the Gateaux second differential of $\mathbb{S}_L^W$ at $0$ with respect to $\langle \cdot,\cdot \rangle_0$,
the linear vector field $Y(\gamma) := - D^2 \mathbb{S}_L^W (0) \gamma$ satisfies
\begin{equation}
\label{Atesi}
D\mathbb{S}_L^W(\gamma) [Y(\gamma)] \leq  - \lambda \|\gamma\|_0^2,
\end{equation}
for every $\gamma$ in a sufficiently small neighborhood of $0$ in $H^1_W([0,1],\R^n)$.
\end{lemma}

\noindent \proof Let $\langle \cdot,\cdot \rangle_0$ be the Hilbert
product
\[
\langle \xi, \eta \rangle_0 := \int_0^1 \frac{\partial^2 L}{\partial v^2} (t,0,0) \xi'(t) \cdot \eta'(t)\, dt + \int_0^1 \xi(t) \cdot \eta(t)\, dt
\]
on $\mathbb{H} := H^1_W([0,1],\R^n)$, and let $\|\cdot\|_0$ be the associated norm. By (\ref{diff2}), the Gateaux Hessian of $\mathbb{S}_L^W$ at $\gamma\in \mathbb{H}$ can be written as
\[
D^2 \mathbb{S}_L^W(\gamma) = A(\gamma) + K(\gamma),
\]
where $A(\gamma)$ and $K(\gamma)$ are the self-adjoint bounded operators on $\mathbb{H}$ defined by
\begin{eqnarray*}
\langle A(\gamma) \xi,\eta \rangle_0 = \int_0^1 \frac{\partial^2 L}{\partial v^2} (t,\gamma(t),\gamma'(t)) \xi'(t) \cdot \eta'(t)\, dt + \int_0^1 \xi(t) \cdot \eta(t)\, dt, \\
\langle K(\gamma) \xi,\eta \rangle_0 = \int_0^1 \Bigl( \frac{\partial^2 L}{\partial q^2} (t,\gamma(t),\gamma'(t)) \xi(t) \cdot \eta(t) + \frac{\partial^2 L}{\partial v \partial q} (t,\gamma(t),\gamma'(t)) \xi(t) \cdot \eta'(t) \\ + \frac{\partial^2 L}{\partial q \partial v} (t,\gamma(t),\gamma'(t)) \xi'(t) \cdot \eta(t) \Bigr)\, dt  - \int_0^1 \xi(t) \cdot \eta(t)\, dt,
\end{eqnarray*}
By (L2'), the operator $A(\gamma)$ is positive and its spectrum is bounded away from zero, uniformly for every $\gamma\in \mathbb{H}$. By our choice of the Hilbert product, $A(0)=I$, so $Y(\gamma) = -\gamma - K(0) \gamma$.

\medskip

\noindent{\em Claim 1. If $(\gamma_h)$ converges to $\gamma$ in $\mathbb{H}$, then $A(\gamma_h)$ converges strongly to $A(\gamma)$.}

\medskip

The sequence $(A(\gamma_h))$ converges strongly to $A(\gamma)$ if and only if for every $\xi\in \mathbb{H}$ and every bounded sequence $(\eta_h)\subset \mathbb{H}$, the sequence
\[
\langle (A(\gamma_h)-A(\gamma)) \xi, \eta_h \rangle_0 = \int_0^1 \left( \frac{\partial^2 L}{\partial v^2} (t,\gamma_h,\gamma_h') - \frac{\partial^2 L}{\partial v^2} (t,\gamma,\gamma') \right) \xi' \cdot \eta_h'\, dt
\]
is infinitesimal (this follows from the fact that a sequence $(\xi_h)$ converges to zero in $\mathbb{H}$ if and only if $\langle \xi_h,\eta_h \rangle_0$ converges to zero for every bounded sequence $(\eta_h)\subset \mathbb{H}$).  By the Cauchy-Schwarz inequality, we have to show that the sequence
\[
 \left( \frac{\partial^2 L}{\partial v^2} (\cdot ,\gamma_h,\gamma_h') - \frac{\partial^2 L}{\partial v^2} (\cdot ,\gamma,\gamma') \right) \xi'
\]
converges to zero in $L^2$. In order to check the latter fact, we may assume that $\gamma_h'$ converges to $\gamma'$ almost everywhere, because of the standard argument involving subsequences (see footnote \ref{subseq}). Then the conclusion follows from the assumption (L1') together with the dominated convergence theorem.

\medskip

\noindent{\em Claim 2. For every $\gamma\in \mathbb{H}$ the operator $K(\gamma)$ is compact, and the map $\gamma \mapsto K(\gamma)$ is continuous with respect to the norm topology on the algebra of bounded operators.}

\medskip

The operator $K(\gamma)$ is compact, for every $\gamma\in \mathbb{H}$, because the single terms appearing in the expression for the associated bilinear form are continuous in $L^2 \times L^2$, in $L^2 \times H^1$, or in
$H^1 \times L^2$ (as in the proof of statement (iii) of Proposition
\ref{reg}). The sequence of self-adjoint operators $(K(\gamma_h))$ converges to $K(\gamma)$ in the norm topology if and only if for every bounded sequence $(\xi_h)\subset \mathbb{H}$ the sequence
\begin{equation}
\label{dadim}
\begin{split}
\langle (K(\gamma_h) - K(\gamma)) \xi_h ,\xi_h \rangle_0 = \int_0^1 \left( \frac{\partial^2 L}{\partial q^2} (t,\gamma_h,\gamma_h') -  \frac{\partial^2 L}{\partial q^2} (t,\gamma,\gamma') \right) \xi_h \cdot \xi_h \, dt \\ + 2 \int_0^1 \left( \frac{\partial^2 L}{\partial v \partial q} (t,\gamma_h,\gamma_h') -  \frac{\partial^2 L}{\partial v \partial q} (t,\gamma,\gamma') \right) \xi_h \cdot \xi_h' \, dt \end{split}
\end{equation}
is infinitesimal (this follows from the fact that the norm $\|K\|_0$ of a self-adjoint operator $K$ is the supremum of $|\langle K \xi,\xi \rangle_0|$ for $\xi$ in the unit ball). Since $(\xi_h)$ is bounded in $L^{\infty}$, in order to prove that the first integral in (\ref{dadim}) tends to zero for $h\rightarrow \infty$, it is enough to show that the sequence
\begin{equation}
\label{l1}
\left( \frac{\partial^2 L}{\partial q^2} (t,\gamma_h,\gamma_h') -  \frac{\partial^2 L}{\partial q^2} (t,\gamma,\gamma') \right)
\end{equation}
converges to zero in $L^1$. Since $(\gamma_h')$ converges to $\gamma'$ in $L^2$, there exists an $L^2$ function $f$ such that $|\gamma_h'| \leq f$ almost everywhere, for every $h\in \N$. Then the dominated convergence theorem and the assumption (L1') imply that the sequence (\ref{l1}) converges to zero in $L^1$.
By the Cauchy-Schwarz inequality, in order to prove that the second integral in (\ref{dadim}) is infinitesimal, it is enough to show that the sequence
\[
\left( \frac{\partial^2 L}{\partial v \partial q} (t,\gamma_h,\gamma_h') -  \frac{\partial^2 L}{\partial v \partial q} (t,\gamma,\gamma') \right) \xi_h
\]
converges to zero in $L^2$. As before, this follows from (L1') by the dominated convergence theorem. This concludes the proof of Claim 2.

\medskip

The gradient vector field of $\mathbb{S}_L^Q$ with respect to the Hilbert
product $\langle \cdot , \cdot \rangle_0$ can be written as
\[
\nabla \mathbb{S}_L^Q (\gamma) = \int_0^1 \frac{d}{ds} \nabla \mathbb{S}_L^Q(s\gamma) \, ds = \int_0^1 D^2 \mathbb{S}_L^Q (s\gamma) \gamma \, ds = \int_0^1 (A(s\gamma) + K(s\gamma)) \gamma \, ds.
\]
Therefore,
\[
D\mathbb{S}_L^Q(\gamma) [-Y(\gamma)] = \langle \nabla S_L^Q(\gamma) , \gamma + K(0) \gamma \rangle_0 = \int_0^1 \langle (A(s\gamma) + K(s\gamma)) \gamma, \gamma + K(0) \gamma \rangle_0 \, ds.
\]
Since the self-adjoint operator $A(\gamma)$ is positive, we can consider its square root $A(\gamma)^{1/2}$ and write
\begin{equation}
\label{Apri} \begin{split}
D\mathbb{S}_L^Q(\gamma) [-Y(\gamma)] = \int_0^1 \langle (A(s\gamma)^{1/2} + K(s\gamma))^2 \gamma, \gamma \rangle_0 \, ds  + \int_0^1 \langle ( A(s\gamma) K(0) \\
+ K(s\gamma) + K(s\gamma) K(0)
- A(s\gamma)^{1/2} K(s\gamma) - K(s\gamma) A(s\gamma)^{1/2} - K(s\gamma)^2) \gamma, \gamma\rangle_0 \, ds. \end{split}
\end{equation}
By Claim 2, $K(\gamma)$ converges to $K(0)$ in norm for $\gamma\rightarrow 0$. Together with Claim 1 and the property of strong convergence stated above, $A(\gamma) K(0) \rightarrow A(0) K(0) = K(0)$ and $A(\gamma)^{1/2} K(\gamma) \rightarrow A(0)^{1/2} K(0) = K(0)$ in norm  for $\gamma\rightarrow 0$. Therefore, the same is true for the operator $K(\gamma) A(\gamma)^{1/2} = ( A(\gamma)^{1/2} K(\gamma) )^*$. We deduce that the norm of
\[
A(s\gamma) K(0)
+ K(s\gamma) + K(s\gamma) K(0)
- A(s\gamma)^{1/2} K(s\gamma) - K(s\gamma) A(s\gamma)^{1/2} - K(s\gamma)^2
\]
converges to zero for $\gamma\rightarrow 0$, uniformly with respect to $s\in [0,1]$. Therefore (\ref{Apri}) implies that
\begin{equation}
\label{Aineq}
D\mathbb{S}_L^Q(\gamma) [-Y(\gamma)] = \int_0^1 \langle (A(s\gamma)^{1/2} + K(s\gamma))^2 \gamma, \gamma \rangle_0 \, ds + o(\|\gamma\|^2_0), \quad \mbox{for } \gamma \rightarrow 0.
\end{equation}

\medskip

\noindent{\em Claim 3. There exists a positive number $\mu$ and a neighborhood of $0$ in $\mathbb{H}$ such that the spectrum of the self-adjoint operator $A(\gamma)^{1/2} + K(\gamma)$ is disjoint from $[-\mu,\mu]$, for every $\gamma$ in such a neighborhood.}

\medskip

Assuming that the contrary is true, we can find sequences $\gamma_h \rightarrow 0$ in $\mathbb{H}$ and $\mu_h \rightarrow 0$ in $\R$ such that $A(\gamma_h)^{1/2} + K(\gamma_h) - \mu_h I$ is not invertible. Equivalently, the self-adjoint operator
\begin{equation}
\label{Aseq}
I + A(\gamma_h)^{-1/2} K(\gamma_h) - \mu_h A(\gamma_h)^{-1/2}
\end{equation}
is not invertible, for every $h\in \N$. Since $A(\gamma_h)$ converges strongly to $A(0)=I$, and since, by (L2'), the spectrum of each of the operators $A(\gamma_h)$ is uniformly bounded away from zero, also $A(\gamma_h)^{-1/2}$ converges strongly to $I$, by the strong continuity of the functional calculus on the space of self adjoint operators. Then, the fact that the operators $K(\gamma_h)$ are compact and converge to $K(0)$ in norm implies that $A(\gamma_h)^{-1/2} K(\gamma_h)$ converges to $K(0)$ in norm. Together with the fact that the sequence $A(\gamma_h)^{-1/2}$ is bounded - again a consequence of the fact that the spectrum of every operator $A(\gamma_h)$ is uniformly bounded away from zero - this implies that the sequence of operators (\ref{Aseq}) converges to $I+K(0)$ in norm. But $I+K(0)=A(0)+K(0)=D^2 \mathbb{S}_L^Q(0)$ is invertible because $0$ is a non-degenerate critical point of $\mathbb{S}_L^Q$. Therefore, the fact that the group of invertible operators is open in the norm-topology implies that the operators (\ref{Aseq}) are invertible for $h$ large enough. This contradiction proves the claim.

\medskip

\noindent{\em Conclusion.} By Claim 3,
\[
\langle ( A(\gamma)^{1/2} + K(\gamma) )^2 \gamma, \gamma \rangle_0 \geq \mu^2 \|\gamma\|^2_0,
\]
for every $\gamma$ in a neighborhood of $0$ in $\mathbb{H}$. Then (\ref{Aineq}) implies that if $0<\lambda<\mu^2$, the inequality (\ref{Atesi}) holds
for every $\gamma$ in a neighborhood of $0$ in $\mathbb{H}$.
\qed

\medskip

\noindent \proof [of Theorem \ref{lyap}] Fix some  $\gamma_0$ in
$H^1_Q([0,1],M)$ which is not a critical point of $\mathbb{S}_L^Q$.
By using a local chart and a smooth cut-off function, it is easy to
construct a smooth vector field $Y_{\gamma_0}$ on $H^1_Q([0,1],M)$
such that $Y_{\gamma_0} (\gamma_0) = - \grad \mathbb{S}_L^Q
(\gamma_0)$. By the continuity of $Y_{\gamma_0}$ and of $D
\mathbb{S}_L^Q$, $\gamma_0$ has an open neighborhood
$\mathscr{U}_{\gamma_0}$ in $H^1_Q([0,1],M)\setminus \mathrm{crit}
(\mathbb{S}_L^Q)$ such that
\begin{equation}
\label{gen}
D\mathbb{S}_L^Q(\gamma) [Y_{\gamma_0}(\gamma)] \leq - \frac{1}{2} \|\grad \mathbb{S}_L^Q (\gamma) \|^2, \quad \forall \gamma\in \mathscr{U}_{\gamma_0}.
\end{equation}
If $\gamma_0$ is a critical point of $\mathbb{S}_L^Q$, it is a smooth curve, so by using some $(U,\varphi)\in \mathscr{F}_Q$ such that
\[
\varphi(t,0) = \gamma_0(t), \quad \forall t\in [0,1],
\]
we can identify the tangent space of $H^1_Q([0,1],M)$ at $\gamma_0$ with the Hilbert space $\mathbb{H} = H^1_W([0,1],\R^n)$, and
a neighborhood of $\gamma_0$ in $H^1_Q([0,1],M)$ with a neighborhood of $0$ in  $\mathbb{H}$.
By using such an identification and by Lemma \ref{loc}, we can find a
Hilbert product $\langle \cdot, \cdot \rangle_0$ on $T_{\gamma_0} H^1_Q([0,1],M)$ and a positive number $\lambda(\gamma_0)$ such that the linear vector field
\[
Y_{\gamma_0} (\gamma):= - D^2 \mathbb{S}_L^Q(\gamma_0) \gamma
\]
associated to the $\langle\cdot,\cdot \rangle_0$-Hessian $D^2
\mathbb{S}_L^Q(\gamma_0)$ of $\mathbb{S}_L^Q$ at $\gamma_0$ satisfies
\[
D\mathbb{S}_L^Q (\gamma) [Y_{\gamma_0}(\gamma)] \leq - \lambda(\gamma_0) \|\gamma\|_0^2,
\]
for every $\gamma$ in an open neighborhood $\mathscr{U}_{\gamma_0}$ of $\gamma_0$ in $H^1_W([0,1],M)$. Since $\grad \mathbb{S}_L^Q$ is locally Lipschitz
(by Proposition \ref{reg} (i)) and since the norms $\|\cdot\|$ and $\|\cdot\|_0$ are locally equivalent, up to the choice of a smaller $\mathscr{U}_{\gamma_0}$ we may assume that
\begin{equation}
\label{crt}
D\mathbb{S}_L^Q(\gamma) [Y_{\gamma_0}(\gamma)] \leq - \mu(\gamma_0) \|\grad \mathbb{S}_L^Q (\gamma) \|^2, \quad \forall \gamma \in \mathscr{U}_{\gamma_0},
\end{equation}
for some positive number $\mu(\gamma_0)$. We may also assume that the open sets $\mathscr{U}_{\gamma_0}$ for $\gamma_0\in \crit \mathbb{S}_L^Q$ are pair-wise disjoint.

By the paracompactness of $H^1_Q([0,1],M)$, there exists a locally finite refinement $\set{\mathscr{V}_j}{j\in \mathscr{J}}$ of the open covering $\set{\mathscr{U}_{\gamma} }{\gamma \in H^1_Q([0,1],M)}$. Let $\Gamma: \mathscr{J} \rightarrow H^1_Q([0,1],M)$ be a function such that $\mathscr{V}_j \subset \mathscr{U}_{\Gamma(j)}$, for every $j\in \mathscr{J}$. Let $\{\chi_j\}_{j\in \mathscr{J}}$ be a smooth partition of unity subordinated to the open covering
$\set{\mathscr{V}_j}{j\in \mathscr{J}}$ (see e.g.\ \cite{lan99} for the
existence of smooth partitions of unity on Hilbert manifolds), and set
\[
Y(\gamma) := \sum_{j\in \mathscr{J}} \chi_j(\gamma) Y_{\Gamma(j)} (\gamma), \quad \forall \gamma\in H^1_Q([0,1],M).
\]
By construction, $Y$ is smooth.
The inequalities (\ref{gen}) and (\ref{crt}) imply that
\begin{equation}
\label{glob}
D\mathbb{S}_L^Q(\gamma)[Y(\gamma)] \leq - \nu(\gamma) \|\grad \mathbb{S}_L^Q (\gamma) \|^2 , \quad \forall \gamma\in H^1_Q([0,1],M),
\end{equation}
where $\nu(\gamma)$ is the minimum among $1/2$ and $\mu(\gamma_0)$ if $\gamma$ belongs to $\mathscr{U}_{\gamma_0}$ for some $\gamma_0 \in \mathrm{crit} \mathbb{S}_L^Q$, $1/2$ otherwise.

We can make $Y$ bounded with respect to the Riemannian metric $\langle \cdot, \cdot \rangle$ by multiplication by a suitable conformal factor: Given a smooth positive function $\chi$ on $[0,+\infty[$ such that $\chi(s)=1$ for $s\leq 1$, and $s\chi(s)$ is bounded and bounded away from zero on $[1,+\infty[$, we set
\[
X(\gamma) := \chi(\|Y(\gamma)\|) Y(\gamma), \quad \forall \gamma\in H^1_Q([0,1],M).
\]
Then $X$ is smooth and bounded, so the completeness of the Riemannian metric of $H^1_Q([0,1],M)$ implies (i). By (\ref{glob}), $\mathbb{S}_L^Q$ is a Lyapunov function for $Y$, hence also for $X$, proving (ii).
Moreover, in a neighborhood of each critical point
$\gamma_0$ the vector field $X$ coincides with the vector field $Y$, which coincides with the linear vector field $\gamma \mapsto -D^2 \mathbb{S}_L^Q(\gamma_0) \gamma$. These facts imply (iii).
The fact that $D\mathbb{S}_L^Q$ is bounded on each sublevel of $\mathbb{S}_L^Q$ (see Remark \ref{bdgrad}) and the properties of $\chi$ imply that every Palais-Smale sequence for $(X,\mathbb{S}_L^Q)$ is also a Palais-Smale sequence for $(Y,\mathbb{S}_L^Q)$.
Finally, since $\mathbb{S}_L^Q$ has finitely critical points in every sublevel, the function $\nu$ is bounded away from zero in every sublevel, so (\ref{glob}) shows that
a Palais-Smale sequence for $(\mathbb{S}_L^Q,Y)$ is also a Palais-Smale sequence for $(\mathbb{S}_L^Q,-\grad \mathbb{S}_L^Q)$. Therefore, (iv) follows
from Proposition \ref{ps}. Since $X$ is smooth, property (v) can be achieved by a suitable generic perturbation, as explained in \cite{ama06m} (see also Remark \ref{morreg}).
\qed
\begin{corollary}
\label{coro}
Under the assumptions of Theorem \ref{lyap}, the Morse complex of $\mathbb{S}_L^Q$ is well-defined up to isomorphism, and its homology is isomorphic to the singular homology of the path space
\[
P_Q([0,1],M) = \set{ \gamma\in C^0([0,1],M) }{(\gamma(0),\gamma(1)) \in Q}.
\]
In particular, problem (\ref{Aele}-\ref{Abdry1}-\ref{Abdry2}) has at least as many solutions of Morse index $k$ as the rank of the singular homology group $H_k(P_Q([0,1],M))$.
\end{corollary}
\noindent\proof Everything follows from the results of Section
\ref{morcom}, together with Theorem \ref{lyap} and the fact that the
inclusion $H^1_Q([0,1],M) \hookrightarrow P_Q([0,1],M)$ is a
homotopy equivalence. \qed
\begin{remark}
{\rm If we drop condition (L0), multiplicity results for the
solutions of (\ref{Aele}-\ref{Abdry1}-\ref{Abdry2}) in terms of the
cup-length of $P_Q([0,1],M)$ have been obtained for the more general
class of Tonelli Lagrangian functions in \cite{af07}. By using the
methods of \cite{af07} and the results of this paper, one could show
that the multiplicity result of Corollary \ref{coro} holds for
Tonelli Lagrangian functions satisfying (L0).}
\end{remark}

\end{document}